\documentclass{article}
\usepackage{amssymb}
\usepackage{amsmath}
\usepackage{epsfig,psfrag,verbatim,amsfonts}
\usepackage[tight]{subfigure}
\begin{document}
\newtheorem{theorem}{Theorem}
\newtheorem{lemma}{Lemma}
\newtheorem{corollary}{Corollary}
\newtheorem{proposition}{Proposition}
\def\squarebox#1{\hbox to #1{\hfill\vbox to #1{\vfill}}}
\newcommand{\qed}{\hspace*{\fill}
\vbox{\hrule\hbox{\vrule\squarebox{.667em}\vrule}\hrule}\smallskip}
\newcommand{\Ray}{{\mbox{\tt Ray}}}
\newcommand{\IGWR}{{\mbox{\tt IGWR}}}
\newcommand{\IGW}{{\mbox{\tt IGW}}}
\newcommand{\GW}{{\mbox{\tt GW}}}
\newcommand{\CC}{\mbox{${\mathcal C}$}} 
\newcommand{\T}{{\mathcal{T}}}
\renewcommand{\v}{{\tt v}}
\newcommand{\BBB}{\mbox{${\mathbb{B}}$}}
\title{{\bf A Central Limit Theorem
for biased random walks 
on Galton-Watson trees}}

\author{Yuval Peres\thanks{
Dept. of Mathematics and Dept. of Statistics, University of California, 
Berkeley. Research  partially supported by MSRI and by NSF grants \#DMS-0104073 and
\#DMS-0244479  }\and
Ofer Zeitouni\thanks{Department of Mathematics, University of Minnesota,
and Depts. of Mathematics and of Electrical Engineering, Technion.
Research partially supported by MSRI and
by NSF grants \#DMS-0302230 and DMS-0503775.}}
\date{June 22, 2006}
\maketitle
\begin{abstract}
Let ${\cal T}$ be a rooted Galton-Watson tree with offspring
distribution $\{p_k\}$ that has $p_0=0$, mean $m=\sum kp_k>1$
and exponential tails.
Consider the
$\lambda$-biased random walk $\{X_n\}_{n\geq 0}$ on  ${\cal T}$;
this is the nearest neighbor random walk which, when at a vertex $v$ with $d_v$ 
offspring, moves closer to the root with probability $\lambda/(\lambda+d_v)$, and moves to
each of the offspring with probability $1/(\lambda+d_v)$. It is known that  this walk has an 
 a.s.\ constant speed $\v=\lim_n |X_n|/n$ 
 (where $|X_n|$ is the distance of $X_n$ from the root), with $\v>0$ for $ 0<\lambda<m$ and $\v=0$ for $\lambda \ge m$.
For all $\lambda \le m$ , we  prove a quenched CLT for $|X_n|-n\v$.
(For $\lambda>m$ the walk is positive recurrent, and there is no CLT.)
The most interesting case by far is $\lambda=m$, where the CLT has the following form:
  for almost every ${\cal T}$, the ratio
$|X_{[nt]}|/\sqrt{n}$  converges in law  as $n \to \infty$ to a deterministic
multiple of the absolute value of a Brownian motion. Our approach to this case
is based on an explicit description of an invariant measure for the walk from the point of view
of the particle (previously, such a measure was explicitly known only for $\lambda=1$)
and the construction of appropriate harmonic coordinates.
\end{abstract}


\section{Introduction and statement of results}
\label{sec-introduction}
\setcounter{figure}{0}

Let ${\cal T}$ be a rooted Galton-Watson tree with offspring
distribution $\{p_k\}$. That is, the numbers 
of offspring $d_v$ of  vertices $v\in {\cal T}$ are i.i.d. random variables,
with $P(d_v=k)=p_k$. Throughout this paper, we assume
that $p_0=0$, and that $m:=\sum kp_k>1$.
In particular, ${\cal T}$ is almost surely an infinite tree.
For technical reasons, we also assume the existence of exponential moments,
that is the existence of some $\beta>1$ such that $\sum \beta^k p_k<\infty$.
We let $|v|$ stand for the distance of a vertex $v$ from the root of
${\cal T}$, and let $o$ denote the root of ${\cal T}$.

We are interested in $\lambda$-biased random walks on the tree ${\cal T}$.
These are Markov chains $\{X_n\}_{n\geq 0}$ with $X_0=o$ and
transition probabilities 
$$ P_{\cal T}(X_{n+1}=w|X_n=v)=
\left\{\begin{array}{ll}
{\lambda}/({\lambda+d_v})\,,& \mbox{\rm if $v$ is an offspring
of $w$}\,,\\
{1}/({\lambda+d_v})\,,& \mbox{\rm if $w$ is an offspring
of $v$}\,.
\end{array}
\right.
$$
Let $\GW$ denote the law of Galton-Watson trees.
Lyons~\cite{Ly} showed that
\begin{itemize}
\item If $\lambda>m$, then for $\GW$-almost every 
${\cal T}$, the random walk $\{X_n\}$ is positive recurrent.
\item if $\lambda=m$, then for $\GW$-almost every 
${\cal T}$, the random walk $\{X_n\}$ is null recurrent.
\item 
if $\lambda<m$, then for $\GW$-almost every
${\cal T}$, the random walk $\{X_n\}$
is  transient.
\end{itemize}
In the latter case, $\lambda<m$, it was later shown in \cite{LPP1} and \cite{LPP2} that
  $|X_n|/n\to \v>0$ almost surely, with a deterministic $\v=\v(\lambda)$ (an explicit expression for $\v$ is known only for
$\lambda=1$). 

Our interest in this paper is mainly in the critical case $\lambda=m$.
Then, $|X_n|/n$ converges to $0$ almost surely.
Our main result is the 
following.
\begin{theorem}
\label{theo-GW}
Assume $\lambda=m$. Then,
there exists a deterministic constant $\sigma^2>0$ such that for
$\GW$-almost every ${\cal T}$, the processes 
$\{|X_{\lfloor nt \rfloor}|/\sqrt{\sigma^2 n}\}_{t\geq 0}$
converges in law to the absolute value of a standard Brownian motion.
\end{theorem}
Theorem \ref{theo-GW} is proved in Section \ref{sec-GWproof} by
coupling 
$\lambda$-biased walks on $\GW$ trees to 
$\lambda$-biased walks
on  auxiliary trees, which have a marked ray emanating
from the root. The ergodic theory of walks on such trees turns out
(in the special case of $\lambda=m$) to be particularly nice.
We develop this model and state the  Central Limit Theorem (CLT) for it,
Theorem \ref{theo-1},
in Section \ref{sec-IGWR}. The proof of 
Theorem \ref{theo-1}, which is based on constructing appropriate martingales
and controlling the associated corrector,
is developed in 
Sections \ref{sec-MSMPT}, \ref{sec-2} and \ref{sec-aux}.

We conclude by noting that when 
$\lambda>m$, the biased random walk is positive recurrent, and 
no CLT limit is possible. On the other hand, \cite{LPP2} proved 
that when $\lambda<m$ and the walk is transient, there exists
a sequence of stationary
regeneration times. Analyzing these regeneration
times, one deduces a quenched invariance principle 
with a proper deterministic
centering, see Theorem \ref{theo-GWtransient}
in Section \ref{sec-remarks} for the statement.
We note in passing that this improves the annealed invariance principle
derived in \cite{piau} for $\lambda=1$.

%

\section{A CLT for trees with a marked ray}
\label{sec-IGWR}
\setcounter{figure}{0}

We consider 
infinite trees $\T$ with
one (semi)-infinite directed path, denoted $\Ray$, 
starting from a distinguished
vertex, called the {\bf root} and denoted $o$.
For vertices 
$v,w\in\T$, we let $d(v,w)$ denote the length of the
(unique) geodesic connecting $v$ and $w$ (we consider the geodesic as 
containing
both $v$ and $w$, and its length as the number of vertices in it
minus one).
A vertex $w$ is an offspring of a vertex $v$ if $d(v,w)=1$ and either
$d(w,\Ray)>d(v,\Ray)$ or $v,w\in \Ray$ and $d(v,o)>d(w,o)$.
In particular, the root is an offspring of its unique neighbor on
$\Ray$.
For any vertex $v\in {\cal T}$,
we let $d_v$ denote the number of offspring of $v$. 


For $v$ a vertex in ${\cal T}$,
let $R_v\in \Ray$ denote the intersection of the geodesic
connecting $v$ to $\Ray$ with $\Ray$, that is
$d(v,R_v)=d(v,\Ray)$.
For $v_1,v_2\in {\cal T}$, let $h(v_1,v_2)$ denote the horocycle
distance between $v_1$ and $v_2$ (possibly negative), 
which is defined
as the unique function $h(v_1,v_2)$ which equals to
$d(x,v_2)-d(x,v_1)$ for all  vertices
$x$ such that both $v_1$ and $v_2$ are descendants of $x$.
(A vertex $w\in\T$ is a descendant
of $v$ if the geodesic connecting
$w$ to $v$ contains an offspring of $v$.)
We also
write
$h(v)=h(o,v)$; The quantity $h(v)$, which may be
either positive or negative, is the {\it level} to which $v$ belongs,
see Figure 1.
\begin{figure}[t]
\psfrag{h=-3}{$h=-3$}
\psfrag{h=-2}{$h=-2$}
\psfrag{h=-1}{$h=-1$}
\psfrag{h=0}{$h=0$}
\psfrag{h=2}{$h=2$}
\psfrag{h=1}{$h=1$}
\psfrag{Ray}{$\Ray$}
\psfrag{zero}{$\pmb o$}
\epsfig{file=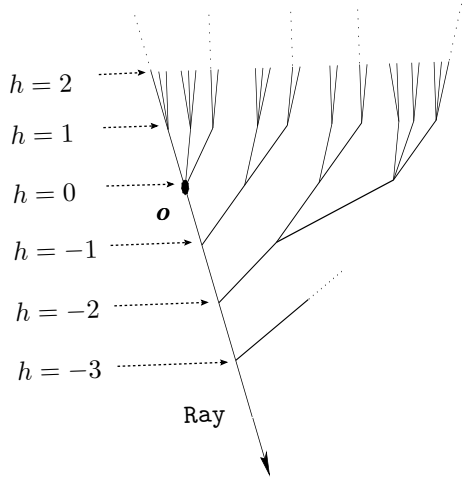,height=6.5cm, width=6.5cm, angle = 2}
\begin{centering}
\caption{Tree, Ray and horocycle distance} 
\end{centering}\label{fig-1}
\end{figure}
Let $D_n(v)$ denote the descendants of $v$ in $\cal T$
at distance $n$ from $v$. Explicitly,
\begin{equation}
	\label{eq-200306b}
	D_n(v)=\{w\in {\cal T}: d(w,v)=h(w)-h(v)=n\}\,.
\end{equation}
We let  $Z_n(v)=|D_n(v)|$ be the number of descendants
of $v$ at level $h(v)+n$. Then   $\{Z_n(v)/m^n\}_{n \ge 1}$ 
forms a martingale and converges a.s., 
as $n\to\infty$, to a random variable denoted
$W_v$. Moreover, $W_v$ has exponential tails, and there 
are good bounds on the rate of convergence, see \cite{At94}.


Motivated by \cite{LPP},
we next describe a  
measure on the collection of trees with marked
rays, which we denote by $\IGW$.
Fix  a vertex $o$ (the root) and a semi-infinite ray, denoted $\Ray$,
emanating from it. Each vertex $v\in \Ray$ with $v\neq o$
is assigned 
independently a size-biased number of offspring, that
is $P_{\IGW}(d_v=k)=kp_k/m$, one of which is
identified with the descendant of $v$ on $\Ray$.
To each offspring of $v\neq o$ not on $\Ray$,
and to $o$,
one attaches an independent Galton-Watson tree of offspring
distribution $\{p_k\}_{k\geq 1}$.
The resulting random
tree $\cal T$ is distributed according to $\IGW$. An alternative 
characterization of $\IGW$ is obtained as follows, see
\cite{LPP} for a similar construction. 
\begin{lemma}
\label{lem-IGW}
Consider the measure $Q_n$
on rooted trees with root $r$, obtained from
$\GW$ 
by size-biasing 
with respect to $|D_n(r)|$ (that is, $dQ_n/d\GW=|D_n(r)|/m^n$).
Choose a vertex $o\in D_n(r)$
uniformly, creating a (finite) ray from $o$ to the root of the
original tree, and extend the ray from $r$ to obtain an infinite ray,
creating thus a random rooted tree with marked ray emanating from the 
new root $o$. Call $\IGW_n$ the distribution thus obtained. Then,
$\IGW$ is the weak limit of $\IGW_n$.
\end{lemma}

Sometimes, we also need to consider trees where the root has no ancestors. 
Often, these will be distributed according
to the Galton-Watson measure $\GW$.
There is however
another important measure that we will use, 
described in \cite{LPP},
namely the {\bf size-biased}  measure
$\widehat{\GW}$
corresponding to $\GW$. It is defined formally by $d\widehat{\GW}/
d\GW= W_o$.
An alternative construction of $\widehat{\GW}$ is by 
sampling, size-biased, a
particular trunk. 

We let 
$\{X_n\}$ denote the $\lambda$-biased
random walk  on the tree $\cal T$, where $\lambda=m$. 
Explicitly,
 given a tree $\cal T$,
$X_n$ is a Markov process with $X_0=o$ and transition probabilities
$$P_{\cal T}(X_{n+1}=u|X_n=v)=\left\{
\begin{array}{ll}
{\lambda}/({\lambda+d_v})\,,& \mbox{\rm if }\;
1=d(u,v)=h(u,v)\,\\
{1}/({\lambda+d_v})\,,&\mbox{\rm if}\;
1=d(u,v)=h(v,u)\,\\
0\,,&\mbox{\rm else}\,.
\end{array}
\right.
$$ 
That is, the walker moves with probability $\lambda/(\lambda+d_v)$
toward the ancestor of $v$ and with probability $1/(\lambda+d_v)$
toward any of the offspring
of $v$.
We recall that the model of $\lambda$-biased random 
walk on a rooted tree is reversible, and possesses
an electric network interpretation, where the conductance
between $v\in D_n(o)$ and an offspring $w\in D_{n+1}(o)$ of $v$
is
$\lambda^{-{n}}$ (see e.g. \cite{LP} for this
representation, and \cite{doyle-snell} for general
background on reversible random walks interpreted in 
electric networks terms). 
With a slight abuse of notation, we 
let $P_{\cal T}^v$ denote the
law, conditional on the given tree $\cal T$ and $X_0=v$,
 on the path $\{X_n\}$. We refer to this law as the
{\it quenched} law. 
Our main result for the $\IGW$ trees is the following.
\begin{theorem}
\label{theo-1}
Under  $\IGW$, the horocycle distance satisfies
a quenched invariance principle. 
That is, for some deterministic $\sigma^2>0$
(see (\ref{eq-121005ss}) below for the value of $\sigma$), for
$\IGW$-a.e. $\cal T$,  the
processes $\{h(X_{\lfloor nt\rfloor
})/\sqrt{\sigma^2 n}\}_{t\geq 0}$  converge 
in distribution to a standard Brownian motion.
\end{theorem}
\section{Martingales, stationary measures, and proof
of Theorem \ref{theo-1}}
\label{sec-MSMPT}
The proof of Theorem \ref{theo-1} takes the bulk of this paper. We describe here
the main steps.
\begin{itemize}
	\item In a first step, we construct in this section 
a martingale $M_t$, whose increments consist of the normalized population size
$W_{X_{t+1}}$ when $h(X_{t+1})-h(X_t)=1$ and 
$-W_{X_{t}}$ otherwise. (Thus, the increments
of the martingale depend on the ``environment as seen from the particle'').
This martingale provides  ``harmonic coordinates'' for the random walk,
in the spirit of \cite{kozlov} and, more recently, \cite{SS} and \cite{BB}.
\item In the next step, we prove an invariance principle for the martingale
	$M_t$. This involves proving a law of large numbers for the associated
	quadratic variation. It is at this step that it turns out that
	$\IGW$ is not so convenient to work with, since the environment viewed
	from the point of view of the particle is not stationary under $\IGW$.
	We thus construct a small modification of $\IGW$, called $\IGWR$,
	which is a reversing measure for the environment viewed from
	the point of view of the particle, and is absolutely continuous 
	with respect to $\IGW$ (see Lemma \ref{lem-rever}). 
	This step uses crucially that $\lambda=m$. Equipped with the measure
	$\IGWR$, it is then easy to prove an invariance principle for
	$M_t$, see Corollary 
	\ref{lem-evp}.
\item In the final step, we introduce the corrector $Z_t$, which is the 
	difference  between 
a constant multiple $1/\eta$ of the harmonic coordinates $M_t$ and
the position of the random walk, $X_t$.  As in \cite{BB}, we seek to show
that the corrector is small, see Proposition \ref{lem-1}. The proof of
Proposition \ref{lem-1} is postponed to Section \ref{sec-2}, and is based
on estimating the time spent by the random walk at any given level.
\end{itemize}

In the sequel (except in Section \ref{sec-GWproof}), 
we often use the letters
$s,t$ to denote time, reserving
the letter $n$ to denote distances on the tree $\cal T$.
Set $M_0=0$ and, if $X_t=v$
for
a vertex $v$ with parent $u$ and offspring  
$Y_1,\ldots,Y_{d_v}$, set
$$
M_{t+1}-M_t=
\left\{\begin{array}{ll}
-W_v, & X_{t+1}=u\\
W_{Y_j},& X_{t+1}=Y_j\,.
\end{array}
\right.
$$
Quenched (i.e., given the realization of the tree), $M_t$ is 
a martingale with respect to the natural filtration
${\cal F}_t=\sigma(X_1,\ldots,X_{t})$, 
as can be seen by using the relation
$W_v=\sum_{j=1}^{d_v} W_{Y_j}/m$. Also, for $v\in \mathcal{T}$,
let $g_v$ denote the geodesic connecting $v$ with $\Ray$ (which by 
definition contains both $v$ and $R_v$), and
set 
$$S_v=\left\{\begin{array}{ll}
\sum_{u\in g_v, u\neq o} W_u, & \mbox{\rm if} \, R_v=o\,,\\
\sum_{ u\in g_v, u\neq R_v} W_u-
\sum_{u\in \Ray,0 \geq h(u)> h(R_v)} W_u, &\mbox{\rm if}
\, R_v\neq o\,.
\end{array}
\right.
$$
Then, $M_t=S_{X_t}$.

Set $\eta=E_{GW}W_o^2 (= E_{\widehat{ GW}} 
W_o)$ and $Z_t=M_t/\eta-h(X_t)$.
Fix 
\begin{equation}
\label{eq-150206}
\alpha=1/3, \epsilon_0<1/100\,,
\delta\in (1/2+\alpha+4\epsilon_0,1-4\epsilon_0).
\end{equation}
(The reason for the
particular choice of constants here will become clearer
in the course of the proof.)
For any integer $t$, let $\tau_t$ denote an integer valued 
random variable,
independent of $\cal T$ and $\{X_s\}_{s\geq 0}$,
uniformly chosen in $[t,t+\lfloor t^\delta\rfloor]$.
We prove in Section \ref{sec-2} the following estimate,
which shows that $M_t/\eta$ is close to $h(X_t)$. The variable $\tau_t$
is introduced here for technical reasons as a smoothing device, that
allows us to consider occupation measures instead of pointwise
in time estimates on probabilities.
\begin{proposition}
\label{lem-1}
With the above notation, for any $\epsilon<\epsilon_0$, 
\begin{equation}
\label{eq-121005kk}
 \lim_{t\to\infty}
P_{\cal T}^o(|Z_{\tau_t}|\geq \epsilon \sqrt{t})
=0\,,\quad
\IGW - a.s.
\end{equation}
Further,
\begin{equation}
\label{eq-100505g}
\lim_{t\to\infty}
P_{\cal T}^o(
\sup_{r,s\leq t, |r-s|<t^\delta}
 | h(X_r)-h(X_s)|>t^{1/2-\epsilon})=0\,,
\; \IGW-a.s.
\end{equation}
\end{proposition}

The interest in the martingale $M_t$ is that we can prove for it a full
invariance principle.
Toward this end, one needs to verify that the
normalized quadratic variation process
\begin{equation}
	\label{eq-200306a}
	V_t=\frac1t \sum_{i=1}^t E_{\cal T}^o
\left( (M_{i+1}-M_i)^2
| {\cal F}_i\right)
\end{equation} converges
$\IGW$-a.s.
Note that
if $X_i=v$ with offspring $Y_1,\ldots,Y_{d_v}$ then
\begin{eqnarray}
\label{eq-var1}
  E_{\cal T}^o\left[ (M_{i+1}-M_i)^2
| {\cal F}_i\right]&=&\frac{m}{m+d_v} W_v^2+
\frac{1}{m+d_{v}}\sum_{j=1}^{d_{v}}W_{Y_j}^2\\
&=&
\frac{1}{m+d_v} \sum_{j=1}^{d_v} W_{Y_j}^2+
\frac{1}{m(m+d_v)}\left(
\sum_{j=1}^{d_v} W_{Y_j}\right)^2=: \mu_v^2\,.
\nonumber
\end{eqnarray}

It turns out that to ensure the convergence of $V_t$, 
it is useful to introduce a new measure on trees, denoted
$\IGWR$, which is absolutely
continuous with respect to the measure $\IGW$, and such that
the ``environment viewed from the point of view of the particle''
becomes stationary under that measure, see Lemma \ref{lem-rever}
below.
The measure $\IGWR$ is similar to $\IGW$, except at the root.
The root $o$
has an infinite path ${v_j}$
 of ancestors, which all possess an independent
number of offspring which is size-biased, that is
$$P(d_{v_j} =k) = kp_k/m\,, \quad \mbox{\rm  for all} \,j,k>0.$$
The number of offspring at the root itself is independent
of the variables just mentioned, and possesses
a distribution which is the
 average of the original and
the size biased
laws, that is:
$$P(d_{o} =k) = (m+k)p_k/(2m)\,, \quad \mbox{\rm for all} \, k>0.$$
All other vertices have the original offspring law.
All these offspring variables are independent. In other words,
$d\IGWR/d\IGW=(m+d_o)/2d_o$. Consequently, we can use
the statements ``$\IGW$-a.s.'' and ``$\IGWR$-a.s.'' interchangeably.

For $v$ a neighbor of $o$, let 
$\theta^v {\cal T}$ denote the tree which is obtained by 
shifting the location
of the root to $v$ and adding or erasing one edge from
$\Ray$ in the only way that leaves an infinite ray emanating 
from the new root. We also write, for an arbitrary vertex $w\in {\cal T}$
with geodesic $g_w=(v_1,v_2,\ldots,v_{|w|-1},w)$ connecting
$o$ to $w$, the shift
$\theta^w{\cal T}=\theta^w\circ\theta^{v_{|w|-1}}\circ \ldots \theta^{v_1}
{\cal T}$. Finally, we set
${\cal T}_t=\theta^{X_t} {\cal T}$. It is evident
that ${\cal T}_t$ is a Markov process, with the location of the random
walk being frozen at the root, and
we write $P_{\cal T}(\cdot)$ for its transition density, that is
$P_{\cal T}(A)=P_{\cal T}({\cal T}_1\in A)$.
What is maybe surprising at first is that 
$\IGWR$ is reversing for this
Markov process. That is, we have.
\begin{lemma}
\label{lem-rever}
The Markov process
${\cal T}_t$ with initial measure $\IGWR$ is stationary and reversible.
\end{lemma}
{\bf Proof of Lemma \ref{lem-rever}}
%
Suppose that ${\cal T}_0$ is picked from IGWR, and
${\cal T}_1$ is obtained from it by doing one step (starting with $X_0=o$) of
the critically biased walk on ${\cal T}_0$, then moving the root to $X_1$ and
adjusting Ray accordingly.
We must show that the ordered pair
$({\cal T}_0, {\cal T}_1)$ has the same law as $({\cal T}_1, {\cal T}_0)$.

Let ${\cal T}_F$ be finite tree of depth $\ell$ rooted at $\rho$,
and let $u,v$ be adjacent internal nodes of ${\cal T}_F$,
at distance $k$ and $k+1$, respectively, from $\rho$ (see figure
2).
\begin{figure}[t]
\psfrag{rho}{$\rho$}
\psfrag{ell}{$\ell$}
\psfrag{u}{$u$}
\psfrag{v}{$v$}
\epsfig{file=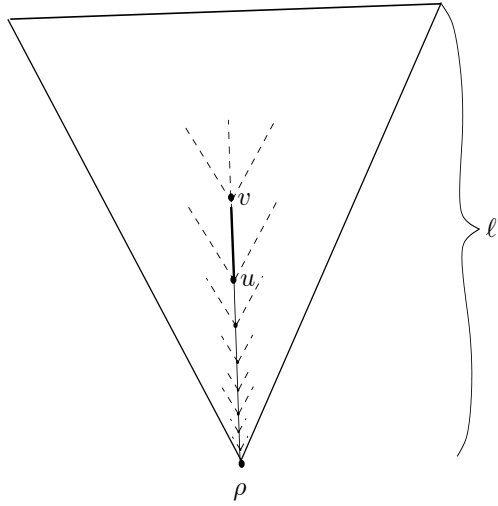,height=6.5cm, width=6.5cm, angle = 2}
\begin{centering}
	\caption{The finite tree ${\cal T}_F$} 
\end{centering}\label{fig-5}
\end{figure}

Let $A({\cal T}_F,u)$ be 
the cylinder set of infinite labeled rooted trees ${\cal T}$ in
the support of $\IGWR$ which locally truncate to ${\cal T}_F$ 
rooted at $u$, that is,
the connected component of the root of ${\cal T}$ among levels between $-k$ and
$\ell-k$ in ${\cal T}$ is identical to ${\cal T}_F$ 
once the root of ${\cal T}$ is identified
with $u$, and $\Ray$\ in ${\cal T}$ goes 
through the vertex identified with $\rho$ in
${\cal T}$. Let
$\{w: \rho \le w <u\}$ denote the set of vertices on the path from
$\rho$ (inclusive) to $u$ (exclusive) in ${\cal T}_F$.
Then
\begin{equation} \label{cyl}
	P_{{\IGWR}}[A({\cal T}_F,u)]= P_{\GW}({\cal T}_F)
\prod_{\{w: \rho \le w <u\}} \Bigl[\frac{d_w}{m}\cdot \frac{1}{d_w} \Bigr]
\frac{m+d_u}{2m} \,,
\end{equation}
where the factors ${d_w}/{m} $ and $(m+d_u)/(2m)$
come from the density of the  $IGWR$ offspring distributions
with respect to the $GW$ offspring distribution,
and the factors $1/d_w$ comes from the uniformity
in the choice of $\Ray$. Thus
\begin{equation} \label{cyl1}
	P_{{\IGWR}}[A({\cal T}_F,u)]= P_{\GW}({\cal T}_F)
	m^{-k-1} (m+d_u)/2 \,,
\end{equation}
and similarly
\begin{equation}
\label{cyl2}
P_{{\IGWR}}[A({\cal T}_F,v)]= P_{\GW}({\cal T}_F) m^{-k-2} (m+d_v)/2 \,.
\end{equation}
Since the transition probabilities for the
critically biased random walk are
$p(u,v)=1/(m+d_u)$ and $p(v,u)=m/(m+d_v)$, we infer from
(\ref{cyl1}) and (\ref{cyl2}) that
$$
P_{{\IGWR}}[A({\cal T}_F,u)] p(u,v) =  P_{{\IGWR}}[A({\cal T}_F,v)] p(v,u)
$$
as required.

\qed

With $V_t$ as in (\ref{eq-200306a}),
the following corollary is of crucial importance.
\begin{corollary}
\label{lem-evp}
\begin{equation}
\label{eq-121005ss}
V_t\to E_{\IGWR} \mu_0^2=:\sigma^2\eta^2\,,\; \IGWR-a.s.
\end{equation}
\end{corollary}
{\bf Proof of Corollary \ref{lem-evp}} 
That $\IGWR$ is absolutely continuous with respect to
$\IGW$ is obvious from the construction.
By Lemma \ref{lem-rever},
$\IGWR$ is invariant and reversible
under the 
Markov dynamics induced by the process ${\cal T}_t$.
Thus,
(\ref{eq-121005ss}) holds  as soon as one checks that
$\mu_0\in L^2(\IGWR)$, which is equivalent to checking
that with $v_i$ denoting the offspring 
of $o$, it holds that $(\sum_{i=1}^{d_o} W_{v_i})^2\in L^1(\IGWR)$.
This in turn is implied by $E_{\GW}(W_o^{2})<\infty$,
which holds  due to \cite{At94}.
\qed


\noindent
{\bf Proof of Theorem \ref{theo-1}}
In what follows, we consider a fixed $\cal T$, with the
understanding that the statements hold true for $\IGW$ almost every 
such tree.
Due to (\ref{eq-121005ss})
and the invariance principle for the Martingale $M_t$,
see \cite[Theorem 14.1]{bil}, it holds that for
$\IGWR$ almost every $\cal T$, 
$\{M_{\lfloor nt\rfloor }/\sqrt{\eta^2
\sigma^2 n}\}_{t\geq 0}$ converges in distribution,
as $n\to\infty$,  to a standard Brownian motion. 
Further, by  \cite[Theorem 14.4]{bil}, 
so does $\{M_{\tau_{nt}}/\sqrt{\eta^2 \sigma^2 n}\}_{t\geq 0}$.
By  (\ref{eq-121005kk}), it then follows
that the finite dimensional distributions
of the process $\{Y^n_t\}_{t\geq 0}=
\{h(X_{\tau_{nt}})/\sqrt{\sigma^2n}\}_{t\geq 0}$
converge, as $n\to\infty$, to those of a standard Brownian
motion. On the other hand, due to 
(\ref{eq-100505g}), the sequence of processes $\{Y_t^n\}_{t\geq 0}$
is tight, and hence converges in distribution to standard
Brownian motion. Applying again \cite[Theorem 14.4]{bil},
we conclude that the sequence of processes
$\{h(X_{\lfloor nt\rfloor })/\sqrt{\sigma^2 n}\}_{t\geq 0}$ converges
in distribution to a standard Brownian motion, as claimed.
\qed


\section{Proof of Proposition \ref{lem-1}}
\label{sec-2}

\noindent
{\bf Proof of  Proposition \ref{lem-1}}
For any tree with root $o$, we write $D_n$ for $D_n(o)$, c.f.
(\ref{eq-200306b}). Recall that $E_{\widehat{\GW}} W_o=\eta$.
For $\epsilon>0$,
let $A_n^{\epsilon}=A_n^\epsilon({\cal T})=
\{v\in D_n: \; |n^{-1} S_v -\eta|>\epsilon\}$, noting that for $\GW$ or
$\widehat{\GW}$ trees, $S_v=\sum_{u\in g_o, u\neq o} W_u$.
We postpone for a moment the proof of the following.
\begin{lemma}
\label{lem-2}
For any $\epsilon>0$
there exists a deterministic $\nu=\nu(\epsilon)>0$ such that
\begin{equation}
\label{eq-oof}
\limsup_{n\to\infty}
\frac1n \log P_{\widehat{ \GW}}\left( 
\frac1n\log\frac
{|A_n^{\epsilon}|}{|D_n|}>-\nu\right)\leq -\nu/2\,,
\end{equation}
and
\begin{equation}
\label{eq-oof1}
\limsup_{n\to\infty}
\frac1n \log P_{\GW}\left( 
\frac1n\log\frac
{|A_n^{\epsilon}|}{|D_n|}>-\nu\right)\leq -\nu/2
\,.
\end{equation}
\end{lemma}
Turning our attention to trees governed
by the measure $\IGW$, for any vertex $w\in \mathcal{T}$ we set
$$S_w^{\Ray}=\sum_{v\in \mathcal{T}\setminus \Ray: v \;\mbox{\rm is
on the geodesic connecting $w$ and $\Ray$}} W_v\,.$$
Let 
$B_n^{\epsilon}(\mathcal{T})=\{w\in \mathcal{T}: d(w,\Ray)=n,
|n^{-1}S_w^{\Ray}-\eta|>\epsilon\}$,
and set
\begin{equation}
\label{eq-190106b}
Q_t(\mathcal{T})=\{w\in \mathcal{T}: d(w,\Ray)
\leq t^{\alpha}\}.\end{equation}
The following proposition will be proved in Section \ref{sec-aux}.
\begin{proposition}
\label{prop-101005a}
\begin{equation}
\label{eq-101005b}
\limsup_{t\to\infty}
P_{\cal T}^o(X_{\tau_t}\in Q_t(\mathcal{T}))=0\,,\;\IGW-a.s.\,.
\end{equation}
\end{proposition}
We can now prove the following.
\begin{lemma}
\label{lem-3}
With the preceding notation, it holds that
for any $\epsilon>0$,
$$\lim_{t\to\infty} P_{\cal T}^o(X_{\tau_t} \in \cup_m 
B_m^{\epsilon}(\mathcal{T})) =0\,,\quad \IGW-a.s.$$
\end{lemma}
{\bf Proof of Lemma \ref{lem-3}}
By  (\ref{eq-101005b}),
\begin{equation}
\label{eq-131005gg}
a_t:=P_{\cal T}^o(X_{\tau_t}\in Q_t(\mathcal{T})
)
\to_{t\to\infty} 0\,,\; \IGW-a.s.
\end{equation}
Letting $\gamma_m^{\epsilon}=\min(t: X_t\in
 B_m^{\epsilon}(\mathcal{T}))$,
we have (using $t+\lceil t^\delta\rceil\leq 2t$),
\begin{equation}
\label{eq-6} P_{\cal T}^o(X_{\tau_t}
 \in \cup_m B_m^{\epsilon}(\mathcal{T})) 
\leq a_t+
\sum_{\ell=t^{\alpha}}^{2t}
P_{\cal T}^o(\gamma_\ell^{\epsilon}\leq 2t)
\,.
\end{equation}
Consider the excursions of $\{X_i\}$ down the $\GW$ trees
whose starting points
are offspring of a vertex in $\Ray$, where an excursion is counted 
between visits to such a starting point.
The event
$\{\gamma_\ell^{\epsilon}\leq 2t\}$ implies
that of the first $2t$ such excursions,
there is at least one excursion that reaches level $\ell-1$ below the
corresponding starting point, at a vertex 
$v$ with $|\ell^{-1} S_v-\eta|>\epsilon$. Therefore, 
with 
$\tau_o=\min\{t>0: X_t=o\}$, for $\ell$ large so that 
$\{x>0: |\ell^{-1} x -\eta|>\epsilon\}
\subset \{x>0: |(\ell-1)^{-1} x-\eta|>\epsilon/2\}$,
\begin{equation}
\label{eq-131005hhhh}
P_{\IGW}^o
(\gamma_\ell^{\epsilon}\leq 2t)\leq
2tP_{\GW}^o(\bar 
\gamma_{\ell-1}^{\epsilon/2}\leq 2t\wedge  \tau_o
)\,,
\end{equation}
where we set
for a $\GW$ rooted  tree,
$\bar \gamma_\ell^{\epsilon/2}=\min\{i>0: X_i\in A_\ell^{\epsilon/2})
\}$.
But, for a $\GW$ rooted tree,
the conductance ${\cal C}({o\leftrightarrow A_{\ell}^{\epsilon/2}})$
from
the root to the vertices 
in $A_\ell^{\epsilon/2}$ is at most $\lambda^{-\ell} |A_\ell^{\epsilon/2}|$.
Note that with $Z_n:=|D_n|m^{-n}$ it holds that $E_\GW(Z_n)=1$ and
$$E_\GW (Z_{n+1}^2)=E_{\GW}(Z_n^2)+
\frac{E_\GW (d_o^2-d_o)}{\lambda^2}(E_\GW(Z_n))^2$$
and hence $E_\GW(Z_\ell^2)\leq c \ell$ for some deterministic
constant $c$. Therefore,
\begin{align*}
P_{\GW}^o(\bar \gamma_{\ell-1}^{\epsilon/2}\leq  \tau_o
)
&\leq 
E_{\GW}(
{\cal C}({o\leftrightarrow A_{\ell-1}^{\epsilon/2}}))\leq 
E_{\GW}(\lambda^{-\ell+1}
|A_{\ell-1}^{\epsilon/2}|)=
E_{\GW}(Z_{\ell-1}
\frac{|A_{\ell-1}^{\epsilon/2}|}{|D_{\ell-1}|})\\
&\leq [E_{\GW}(Z_{\ell-1}^2)]^{1/2} 
[E_{\GW}( (\frac{|A_{\ell-1}^{\epsilon/2}|}{|D_{\ell-1}|})^2))]^{1/2}
\leq
e^{- \nu(\epsilon/2) \ell/4}\,.
\end{align*}
for $\ell$ large, where Lemma \ref{lem-2} was used in the last inequality.
Combined with (\ref{eq-131005hhhh}), we conclude that
$$\sum_{\ell=t^\alpha}^{2t}
P_{\IGW}^o
(\gamma_\ell^{\epsilon}\leq 2t)\leq
e^{- \nu(\epsilon/2)t^\alpha/8}\,.$$
By Markov's inequality and the Borel-Cantelli lemma, this implies that
$$ \limsup_{t\to\infty}
e^{\nu(\epsilon/2) t^\alpha/16}\sum_{\ell=t^\alpha}^{2t}
P_{\cal T}^o
(\gamma_\ell^{\epsilon}\leq 2t)=0\,, \quad \IGW-a.s.
$$
Substituting in (\ref{eq-6}) and using (\ref{eq-131005gg}), 
 one concludes the proof of
Lemma \ref{lem-3}.
\qed

\noindent
{\bf Proof of Lemma \ref{lem-2}} 
Recall the construction of the measures
$\widehat{\GW}$ and $\widehat{\GW}_*$, see \cite[Pg 1128]{LPP}.
Note that $\widehat{\GW}_*$ is a measure on rooted trees with
a marked ray emanating from the root. 
We let $v_n^*$ denote the marked vertex at 
distance $n$ from the root.

By \cite[(2.1),(2.2)]{LPP}, and denoting 
by $\T_n$ the first $n$ generations of
the tree $\T$, it holds that
$$\widehat{\GW}_*(v_n^*\in A_n^{\epsilon})=
E_{\widehat{\GW}}\left(\frac{1}{|D_n|}
\sum_{v\in D_n} P_{\widehat{\GW}}(v\in A_n^{\epsilon}|\T_n)\right)\,.$$
We show below that there exists $\delta_1=\delta_1(\epsilon)>0$ 
such that
\begin{equation}
\label{eq-roscoff1}
\widehat{\GW}_*(v_n^*\in A_n^{\epsilon})\leq e^{-2\delta_1 n}\,.
\end{equation}
We assume that (\ref{eq-roscoff1}) has been proved, and complete
the proof of the lemma. 
By Markov's inequality, (\ref{eq-roscoff1}) implies that
\begin{eqnarray}
\label{eq-141005a}
&& P_{\widehat{\GW}}\left(E_{\widehat{\GW}}
\left(\frac{|A_n^{\epsilon}|}{|D_n|} 
\mid \T_n\right)\geq e^{-\delta_1 n}\right)\\
&=& P_{\widehat{\GW}}\left(\frac{1}{|D_n|} \sum_{v\in D_n}
P_{\widehat{ \GW}}(v\in A_n^{\epsilon}
\mid \T_n)\geq e^{-\delta_1 n}\right)\leq 
\frac{e^{-2\delta_1 n}}{e^{-\delta_1 n}}=e^{-\delta_1 n}\,.
\nonumber \end{eqnarray}
We thus get
\begin{eqnarray*}
&& P_{\widehat{\GW}}\left(\frac{|A_n^{\epsilon}|}{|D_n|} 
>e^{-\delta_1 n/2}\right)=
E_{\widehat{\GW}}\left(P_{\widehat \GW}\left(\frac{|A_n^{\epsilon}|}{|D_n|} 
>e^{-\delta_1 n/2}\mid \T_n\right)\right)\\
&
\leq& 
E_{\widehat{\GW}}\left(E_{\widehat \GW}\left(\frac{|A_n^{\epsilon}|}{|D_n|} 
\mid \T_n\right)
e^{\delta_1 n/2}
\right)\\
&
\leq & e^{-\delta_1 n/2}+e^{\delta_1 n/2}
P_{\widehat{\GW}}\left(E_{\widehat \GW}\left(\frac{|A_n^{\epsilon}|}{|D_n|} 
\mid \T_n\right)\geq e^{-\delta_1 n}\right)
\leq 2e^{-\delta_1 n/2}\,,
\end{eqnarray*}
where Markov's inequality was used in the first inequality
 and (\ref{eq-141005a}) in the last.
This proves (\ref{eq-oof}). While (\ref{eq-oof1}) could be
proved directly, one notes that, with $r>1$ such that
$p_1m^{r-1}<1$,
\begin{eqnarray*}
&&P_{ \GW}\left( 
\frac1n\log\frac
{|A_n^{\epsilon}|}{|D_n|}>-\nu\right)=
E_{\widehat{ \GW}}\left( W_o^{-1} 
{\bf 1}_{\frac1n\log\frac
{|A_n^{\epsilon}|}{|D_n|}>-\nu}\right)
\\
&\leq&
(E_{\widehat \GW}W_o^{-r})^{1/r}
\left(P_{\widehat{ \GW}}\left(  
\frac1n\log\frac
{|A_n^{\epsilon}|}{|D_n|}>-\nu\right)
\right)^{1-1/r}\,,
\end{eqnarray*}
where  H\"{o}lder's inequality with exponent $r>1$ was used.
Since
$E_{\widehat{\GW}}W_o^{-r}=E_{\GW}(W_o^{-(r-1)})<\infty
$ by  \cite[Theorem 1]{NV}, 
(\ref{eq-oof1}) follows from (\ref{eq-oof}).


It remains to prove (\ref{eq-roscoff1}). We use the following: 
Since $$(E_{\widehat{\GW}} e^{\xi W_o})^2\leq
E_{\GW}(W_o^2)
E_{\GW} e^{2\xi W_o}<\infty$$ for some $\xi>0$, where the 
last inequality is due to \cite{At94},
it follows that there exists a $\xi>0$ such that
\begin{equation}
\label{eq-roscoff2}
E_{\widehat{\GW}_*} e^{\xi W_o}= 
E_{\widehat{\GW}} e^{\xi W_o}<\infty. 
\end{equation}
For a marked vertex $v_k^*$, we let $\tilde Z_n^{v_k^*}$ 
denote the size of
the subset of vertices in $D_{n}(v_k^*)$ whose ancestral line does not
contain $v_{k+1}^*$,
and we define $\tilde W_k$ as the a.s. limit (as $n\to\infty$) 
of $\tilde Z_n^{v_k^*}/m^n$, which exists by the standard martingale 
argument. Note that
by construction, for $k<n$, with $W_k=W_{v_k^*}$,
\begin{equation}
	\label{eq-200306c}
	W_k=\tilde W_k+\frac{\tilde W_{k+1}}{m}+
\ldots
+\frac{\tilde W_{n-1}}{m^{n-k-1}}+
\frac{W_n}{m^{n-k}}\,.
\end{equation}
Therefore,
$$S_{v_n^*}=
\sum_{k=0}^{n-1} \tilde W_{k} C_k+ W_n C_n\,,$$
where $C_k=1+1/m+(1/m)^2+\ldots+(1/m)^k$. 
Due to 
(\ref{eq-roscoff2}), we have the existence of a $\delta_2>0$ such
that
\begin{equation}
\label{eq-roscoff3}
P_{\widehat{\GW}_*}(|W_n C_n|>\epsilon n/4)\leq
P_{\widehat{\GW}}(|W_o |>(1-1/m)\epsilon n/4)\leq e^{-\delta_2 n}.
\end{equation}
Also, 
\begin{eqnarray}
\label{eq-roscoff4}
P_{\widehat{\GW}_*}(\sum_{k=0}^{n-1} \tilde W_k[C_\infty-C_k]
>\epsilon n/4)
&=&
P_{\widehat{\GW}_*}(\sum_{k=0}^{n-1} \frac{\tilde W_k}{m^{k+1}(1-1/m)}
>\epsilon n/4)\nonumber\\
&\leq &n
P_{\widehat{\GW}}( \tilde W_o
>c_{\epsilon, m} n )\leq e^{-\delta_2 n}\,,
\end{eqnarray}
for some constant $c_{\epsilon,m}$,
where 
(\ref{eq-roscoff2}) was used in the second inequality.
On the other hand, 
$$ \eta= E_{\widehat{\GW}_*} W_{k}= E_{\widehat{\GW}} [C_\infty 
\tilde{W}_0]\,,$$
where the first equality follows from the construction of $\widehat{\GW}$ and
the definition of $\eta$, and the second from (\ref{eq-200306c}).
The random variables 
$\tilde W_{k}$ are i.i.d. by construction under
$\widehat{\GW}_*$. Therefore, using (\ref{eq-roscoff3}) and
(\ref{eq-roscoff4}), 
$$P_{\widehat{\GW}_*}\left(\left|\frac{S_{v_n^*}}{n}-\eta\right|>\epsilon
\right)\leq 2e^{-\delta_2 n}+
P_{\widehat{\GW}_*}
\left(\frac1n \sum_{k=0}^{n-1} [C_\infty \tilde W_{k}-\eta]>
\frac{\epsilon}{2}\right)\,.$$
Standard large deviations (applied to the sum of i.i.d. random variables
$\tilde W_k$ that possess exponential moments) together with
(\ref{eq-roscoff2}) now yield (\ref{eq-roscoff1}) and complete the
proof of Lemma \ref{lem-2}.\qed

Continuing with the proof 
of Proposition \ref{lem-1},
let $v_n$ denote the vertex on $\Ray$ with $h(v_n)=-n$.
By the same construction 
as in the course
of the proof of Lemma \ref{lem-2}, it holds that
\begin{equation}
\label{eq-131005tt}
S_{v_n}/n\to_{n\to\infty} -\eta\,, \quad \IGW-a.s..
\end{equation}
Let $R_{t}=R_{X_{\tau_t}}$. Note that 
$S_{X_{\tau_t}}=-S_{R_t}
+S_{X_{\tau_t}}^{\Ray}$. Thus,
$$ |Z_{\tau_t}|
\leq  |S_{R_t}/\eta + |h(R_t)||+ 
|S^{\Ray}_{X_{\tau_t}}/\eta- h(R_t,X_{\tau_t})| \,.$$
Note that since the random walk restricted to $\Ray$
is transient, $h(R_t)\to_{t\to\infty}-\infty$, 
and hence by (\ref{eq-131005tt}), 
$S_{R_t}/\eta  |h(R_t)| \to -1$.
Therefore,
for any positive $\epsilon_1$,
for all large $t$, using that $\tau_t\leq 2t$,
it follows that $|S_{R_t}/\eta + |h(R_t)|| 
\leq \epsilon_1  \sup_{s\leq 2t} |M_s|$.
Similarly,
for any $\epsilon_1<\epsilon$,
on the event $X_{\tau_t}\not\in \cup_m 
B_m^{\epsilon_1}({\cal T})$,  
it holds 
that for large $t$,
$|S^{\Ray}_{X_{\tau_t}}/\eta- h(R_t,X_{\tau_t})|
\leq \sup_{s\leq 2t} |M_s|$ for all $t$ large.
Thus, for such $\epsilon_1$,
$|Z_{\tau_t}|\leq 2\epsilon_1  \sup_{s\leq 2t} |M_s|$ 
for all $t$ large.
From Lemma \ref{lem-3}, 
\begin{equation}
\label{eq-071105aaaa}
\limsup_{t\to\infty}
P^o_{\cal T}(
X_{\tau_t} \in \cup_m B_m^{\epsilon_1}({\cal T}) =0\,.
\end{equation} 
But, since the normalized 
increasing process $V_t$ is $\IGWR$-a.s. bounded,
standard Martingale inequalities imply that
$$ \lim_{\epsilon_1\to 0} \limsup_{t\to\infty} P_{\cal T}^o(
 \sup_{s\leq t} |M_t|>\epsilon \sqrt{t}/2\epsilon_1) =0\,.$$
It follows that
$$\lim_{t\to\infty}
P_{\cal T}^o(|Z_{\tau_t}|\geq \epsilon\sqrt {t})=0\,,
$$
as claimed.

%
%
%
The proof of
(\ref{eq-100505g}) is provided  in Section
\ref{sec-aux}, see (\ref{eq-oof2}).
This completes the proof of Proposition \ref{lem-1}.
\qed
\section{Auxiliary computations and proof of (\ref{eq-100505g})}
\label{sec-aux}

We begin by an a-priori annealed estimate on the displacement
of the random walk in a GW tree.
\begin{lemma}
\label{lem-111005a}
For any $u,t\geq 1$, it holds that   
\begin{equation}
\label{eq-071005hhh}
P_{\GW}^o(|X_i|\geq u\; 
\mbox{\rm for some \,$i\leq t$})
\leq 4t e^{-u^2/2t}\,.
\end{equation}
\end{lemma}
{\bf Proof of Lemma \ref{lem-111005a}}
Throughout, we write $|v|=d(v,o)$.
Let $\mathcal{T}_u$ 
denote the truncation of the tree $\cal T$ at level $u$,
and let $\mathcal{T}^*$ denote the graph obtained from $\T_u$
by adding an extra vertex (denoted $o^*$) and
connecting it to all vertices in $D_{u}$.
Let $X_s^*$ denote the random walk on $\mathcal{T}^*$, with 
$$
P_{\cal T}(X_{i+1}^*=w|X_i^*= v)=
\left\{\begin{array}{ll}
P_{\cal T}(X_{i+1}=w|X_i=v),&
 \mbox{\rm if}\, v\not\in D_u,\\
1/2,&
 \mbox{\rm if $v\in D_u$\, and $d(v,w)=1$},\\
1/|D_{u}| ,&
 \mbox{\rm if $v=o^*$\, and $d(v,w)=1$}
\,.\end{array}\right.$$
Then,
\begin{eqnarray}
\label{eq-071005h}
&&P_{\GW}^o(|X_i|\geq u\; 
\mbox{\rm for some \,$i\leq t$})
=
P_{\GW}^o(|X_i^*|= u\; 
\mbox{\rm for some \,$i\leq t$})\nonumber\\
& \leq &\sum_{i=1}^t
P_{\GW}^o(|X_i^*|=u) 
\leq 2 \sum_{i=1}^{t+1}
P_{\GW}^o(|X_i^*|=o^*)\,.
\end{eqnarray}
By the Carne-Varopoulos bound, see 
\cite{carne,var},
\cite[Theorem 12.1]{LP},
$$
P_{\cal T}^o(|X_i^*|=o^*)\leq 2 \sqrt{\lambda^{-u}|D_{u}|/d_o}
e^{-u^2/2i}\,.$$
Hence, since $E_{GW}|D_{u}|=\lambda^{u}$,
$$ 2 \sum_{i=1}^{t+1}
P_{\GW}^o(|X_i^*|=o^*)
\leq 4t e^{-u^2/2t}\,.
$$
Combining the last estimate with
(\ref{eq-071005h}), we get (\ref{eq-071005hhh}). \qed

We get the following.
\begin{corollary}
\label{cor-111005a}
It holds that
\begin{equation}
\label{eq-111005hhh}
P_{\IGWR}^o(|h(X_i)|\geq u\; 
\mbox{\rm for some \,$i\leq t$})
\leq 8t^3 e^{-(u-1)^2/2t}\,.
\end{equation}
and
\begin{equation}
\label{eq-131005hhh}
P_{\IGW}^o(|h(X_i)|\geq u\; 
\mbox{\rm for some \,$i\leq t$})
\leq 16t^3 e^{-(u-1)^2/2t}\,.
\end{equation}
\end{corollary}
{\bf Proof of Corollary \ref{cor-111005a}}
We begin by estimating $P_{\IGWR}^o(h(X_i)\geq u)$.
Note that, decomposing according to the
last visit to the level $0$, 
\begin{eqnarray*}
&&P_{\IGWR}^o(h(X_i)\geq u) \\
&\leq& 
P_{\IGWR}^o(\exists j<i:\,
h(X_i)-h(X_j)\geq u\,,
 h(X_t)-h(X_j)>0\, \forall
t\in \{j+1,\ldots,i\})\nonumber \\
&\leq &
\sum_{j=0}^{i-1}P_{\IGWR}^o(
h(X_i)-h(X_j)\geq u\,,
 h(X_t)-h(X_j)>0\, \forall
t\in \{j+1,\ldots,i\})\,.
\end{eqnarray*}
Using the stationarity of $\IGWR$, we thus get
\begin{eqnarray}
\label{eq-160606a}
&&P_{\IGWR}^o(h(X_i)\geq u) \\
&\leq &\sum_{j=0}^{i-1} P_{\IGWR}^o(h(X_{i-j})\geq u, h(X_s)>0
\; \forall s\in \{1,\ldots,i-j\})\,,\nonumber\\
&\leq &
i \max_{r\leq i}
P_{\IGWR}^o(h(X_{r})\geq u, h(X_s)>0
\;\forall  s\in \{1,\ldots,r\})\nonumber.
\end{eqnarray}
On the other hand,
for $r,u>1$, 
\begin{equation}
	\label{eq-200306d}
	P_{\IGWR}^o(h(X_{r})\geq u, h(X_s)>0
\;\forall  s\in \{1,\ldots,r\})
\leq
P_{\GW}^o(h(X_{r})\geq u-1)\,,
\end{equation}
because reaching level $u$ before time $r$  and 
before returning to the root or visiting $\Ray$
requires reaching level $u$ from one of the offspring of the root 
before returning to the root. Substituting in
(\ref{eq-160606a}) we get
\begin{equation}
\label{eq-111005yy}
P_{\IGWR}^o(h(X_i)\geq u)
\leq  i 
\max_{r\leq i} P_{\GW}^o(h(X_{r})\geq u-1)
\leq 4i^2 e^{-(u-1)^2/2i}\,,
\end{equation}  
where (\ref{eq-071005hhh}) was used in the last inequality. 
It follows from the above that
\begin{equation}
\label{eq-111005hhu}
P_{\IGWR}^o(h(X_i)\geq u\; 
\mbox{\rm for some \,$i\leq t$})
\leq 4t^3 e^{-(u-1)^2/2t}\,.
\end{equation}

Recall the process ${\cal T}_s=\theta^{X_s}{\cal T}$, which is reversible
under $P_{\IGWR}$, and note that
$h(X_i)-h(X_0)$ is a measurable function, say
$H$, of $\{{\cal T}_j\}_{0\leq j\leq i}$
(we use here that for $\IGWR$-almost every ${\cal T}$, and vertices 
$v,w\in {\cal T}$, one has $\theta^v{\cal T}\neq \theta^w{\cal T}$. 
Further, with $\widehat{\cal T}_j:={\cal T}_{i-j}$, it holds that
$H(\{\widehat {\cal T}_j\}_{0\leq j\leq i})=
-H(\{ {\cal T}_j\}_{0\leq j\leq i})$. Therefore,
$$P_{\IGWR}^o(h(X_i)\leq -u)=P_{\IGWR}^o(h(X_i)\geq u)\,.$$
Applying 
(\ref{eq-111005yy}), one concludes that
\begin{equation}
\label{eq-111005huu}
P_{\IGWR}^o(h(X_i)\leq -u\; 
\mbox{\rm for some \,$i\leq t$})
\leq 4t^3 e^{-(u-1)^2/2t}\,.
\end{equation}
Together with (\ref{eq-111005hhu}), the proof of 
(\ref{eq-111005hhh}) is complete. To see
(\ref{eq-131005hhh}), note that
$\IGW$ is absolutely continuous with respect  to $\IGWR$, with 
Radon-Nikodym derivative uniformly bounded by $2$.
\qed

We can now give the\\
{\bf Proof of (\ref{eq-100505g})}
The increments $h(X_{i+1})-h(X_i)$ are stationary 
under $P^o_{\IGWR}$. Therefore, by (\ref{eq-111005hhh}),
for any $\epsilon$ and
$r,s\leq t$ with $|r-s|\leq t^\delta$,
$$P_{\IGWR}^o(| h(X_r)-h(X_s)|>t^{1/2-\epsilon})
=
P_{\IGWR}^o(| h(X_{r-s})|>t^{1/2-\epsilon})
\leq 8 t^3 e^{-t^{1-\delta-2\epsilon}}\,.$$
Therefore, by Markov's inequality, for all $t$ large,
$$ P_{\IGWR}\left(P_{\cal T}^o\left(| h(X_{r-s})|>t^{1/2-\epsilon}\right)
\geq t^{-2} e^ {-t^{1-\delta-\epsilon}}\right)\leq 
e^ {-t^{1-\delta-\epsilon}}\,.$$
Consequently,
$$ P_{\IGWR}\left(P_{\cal T}^o\left(
\sup_{r,s\leq t, |r-s|<t^\delta}
| h(X_r)-h(X_s)|>t^{1/2-\epsilon}\right) 
\geq e^ {-t^{1-\delta-\epsilon}}\right)\leq 
 e^ {-t^{1-\delta-\epsilon}}\,.$$
It follows that 
\begin{equation}
\label{eq-oof2}
\limsup_{t\to\infty}
\frac{ P_{\cal T}^o\left(
\sup_{r,s\leq t, |r-s|<t^\delta}
| h(X_r)-h(X_s)|>t^{1/2-\epsilon}\right) }
{ e^ {-t^{1-\delta-\epsilon}}} \leq 1\,,\quad 
\IGWR-a.s.,\end{equation}
completing the proof of (\ref{eq-100505g}) since the measures
$\IGWR$ and $\IGW$ are mutually absolutely continuous. \qed

We next control the expected number of visits to 
$D_n$ during one excursion
from the root of a $\GW$ tree. 
We recall that $T_o=\min\{n\geq 1: X_n=o\}$.
\begin{lemma}
\label{lem-071005b}
Let ${\cal N}_o(n)=\sum_{i=1}^{T_o} {\bf 1}_{X_i\in D_n}$. There
exists a constant $C$ independent of $n$ such that
\begin{equation}
\label{eq-071005c}
E^o_{\GW}({\cal N}_o(n)|d_o)\leq Cd_o\,\quad \mbox{\rm and }\;
E^o_{\widehat{\GW}}({\cal N}_o(n)|d_o)\leq Cd_o\,.
\end{equation}
Further, 
\begin{equation}
\label{eq-071105ff}
\limsup_{n\to\infty}
E^o_{\cal T}({\cal N}_o(n))<\infty\,,
\quad \GW-a.s.
\end{equation}
%
\end{lemma}
{\bf Proof of Lemma \ref{lem-071005b}}
We begin by conditioning on the tree $\cal T$, 
and fix a vertex $v\in D_n$.
Let $\Gamma_v$ denote the number of visits to $v$ before $T_o$. Then,
$$E^o_{\cal T}(\Gamma_v)=P_{\cal T}^o(T_v<T_o)E_{\cal T}^v(\Gamma_v)\,.$$
Note that the walker  performs, on the ray connecting $o$ and $v$, a
biased random walk with holding times. Therefore, by standard computations,
$$P_{\cal T}^o(T_v<T_o)= \frac{1}{d_o[1+\lambda+\lambda^2+
\ldots+\lambda^{n-1}]}\,,$$
and, when starting at $v$, $\Gamma_v$ is  a Geometric random variable
with parameter 
$\lambda^n/[(\lambda+d_v)(1+\lambda+\lambda^2+\ldots+\lambda^{n-1})]$.
Therefore, for some deterministic constant $C$,
$$E^o_{\cal T}(\Gamma_v)\leq C \lambda^{-n}d_v\,.
$$
Thus,
\begin{equation}
\label{eq-071105gg}
E^o_{\cal T}({\cal N}_o(n))\leq C 
\sum_{v\in D_n} \lambda^{-n} d_v \,.
\end{equation}
Since the random variables
$d_v$ are i.i.d., independent of $D_n$,
 and possess exponential moments,
and since $|D_n|\lambda^{-n}\to_{n\to \infty} W_o<\infty$,
it holds that 
$$\limsup_{n\to\infty} 
\sum_{v\in D_n} \lambda^{-n} d_v<\infty \,.$$
Together with  (\ref{eq-071105gg}), this proves
(\ref{eq-071105ff}).
Further, it follows from (\ref{eq-071105gg}) that
$$E^o_{\GW}\left(
{\cal N}_o(n)\,\Big| d_o\right)\leq C\lambda^{-n} E_\GW
\left(|D_n|\,\Big| d_o\right)= Cd_o\,.$$
The proof for $\widehat{\GW}$ is similar.
\qed

We return to $\IGW$ trees. 
Recall that  $Q_t(\mathcal{T})=\{w\in \mathcal{T}: d(w,\Ray)
\leq t^{\alpha})$, and set
$N_t(\alpha)=\sum_{i=1}^t {\bf 1}_{X_i\in Q_t(\mathcal{T})}\,.$
\begin{lemma}
\label{lem-071005aa}
For each $\epsilon>0$ it holds that for all $t$ large enough,
\begin{equation}
\label{eq-071005e}
E^o_{\IGW}(N_t(\alpha))\leq  t^{1/2+\alpha+\epsilon}\,.
\end{equation}
\end{lemma}
{\bf Proof of Lemma \ref{lem-071005aa} }
Let $U_t=\min\{h(X_i):i\leq t\}$ 
and $t_{\epsilon}=\lceil t^{1/2+\epsilon/4}\rceil$. 
By (\ref{eq-131005hhh}), for $t$ large,
\begin{equation}
\label{eq-071005k}
P_{\IGW}^o(U_t\leq - t_{\epsilon})
\leq 16 t^3e^{-t^{\epsilon/2}/3}\,.
\end{equation}
Let $\xi_i=\min\{s: h(X_s)=-i\}$.
It follows from (\ref{eq-071005k}) that for all $t$ large,
\begin{eqnarray}
\label{eq-071005l}
E^o_{\IGW}(N_t(\alpha))&\leq & 
1+
E^o_{\IGW}(N_t(\alpha); U_t> -t_{\epsilon}) \nonumber \\
&\leq &
1+
E^o_{\IGW}(N_t(\alpha); \xi_{t_{\epsilon}}\geq t)\,. 
\end{eqnarray}

For all $k\geq 0$,
let $v_k$ be the unique vertex on $\Ray$ 
satisfying  $h(v_k)=-k$, and set
$d_k=d_{v_k}$.
We next claim that there exists a constant $C_1=C_1(\epsilon)$ independent
of $t$  such that,
with
$$\Upsilon_{t,\epsilon}:=
\{\max_{k\in[0,t_\epsilon]} d_k \leq C_1 (\log t_\epsilon)\}\,,
$$
it holds that
\begin{equation}
	\label{eq-150206ea}
P_{\IGW}(\Upsilon_{t,\epsilon}^c)\leq \frac1t\,.
\end{equation}
Indeed, with  $\beta'=1+(\beta-1)/2>1$,
\begin{eqnarray}
	\label{eq-150206eb}
	P_{\IGW}(\Upsilon_{t,\epsilon}^c)&\leq &t_{\epsilon} P_{\IGW}(d_0>C_1
	\log t_{\epsilon})\nonumber\\
	&\leq&
	\frac{t_{\epsilon}}{m}
	\sum_{j=
	C_1 \log t_{\epsilon}}^\infty jp_j
	\leq 
	\frac{t_{\epsilon} (\beta')^{-C_1 \log t_{\epsilon}}}{m}
	\sum_{j=1
	}^\infty jp_j(\beta')^j\,,
\end{eqnarray}
from which (\ref{eq-150206ea}) follows if $C_1$ is large enough since
$\sum \beta^j p_j<\infty$ by assumption.
Combined with the fact that $N_t(\alpha)\leq t$ and
(\ref{eq-071005l}), we conclude that for such $C_1$,
\begin{equation}
	\label{eq-150206ec}
	E^o_{\IGW}(N_t(\alpha))\leq  
	2+
E^o_{\IGW}(N_t(\alpha); \xi_{t_{\epsilon}}\geq t; \Upsilon_{t,\epsilon})\,. 
\end{equation}

For the next step, let $\theta_0=0$ and, for $\ell\geq 1$,
let $\theta_\ell$ denote the $\ell$-th visit to $\Ray$,
that is
$\theta_\ell=\min\{t>\theta_{\ell-1}: X_t\in \Ray\}$. Let
$H_\ell=X_{\theta_\ell}$ denote the skeleton of $X_i$ on $\Ray$.
Note that $h_\ell=h(H_\ell)$ 
is a (biased) random walk in random environment with holding times;
that is, 
\begin{equation}
	\label{eq-150206ed}
	P(h_{\ell+1}=j|h_\ell=k)=
\left\{\begin{array}{ll}
{\lambda}/({\lambda+d_k})\,,& j=k-1,\\
{1}/({\lambda+d_k})\,,& j=k+1,\\
{(d_k-1)}/({\lambda+d_k})\,,&j=k\,.
\end{array}
\right.
\end{equation}
Let
$h_\ell^*$ denote the homogeneous 
Markov chain on $\mathbb{Z}$ with $h_0^*=0$ and transitions as in
(\ref{eq-150206ed})
corresponding
to a homogeneous environment with $d_k=C_1 \log t_{\epsilon}$, and
set $\eta_i=\min\{\ell:h_\ell=-i\}$ and
$\eta_i^*=\min\{\ell:h_\ell^*=-i\}$. The chain $h_\ell^*$ possesses the same 
drift as the chain $h_\ell$, and on the event 
$\Upsilon_{t,\epsilon}$,
its holding times dominate those of the latter chain. Therefore,
$$
{\bf 1}_{\Upsilon_{t,\epsilon}}
P_{\cal T}^o(\eta_{t_{\epsilon}}>m)\leq 
P(\eta_{t_{\epsilon}}^*>m)
\,.$$
Further, setting $\bar\theta_0=0$ and, for $j\geq 1$,
using
$\bar\theta_j=\min\{i>\bar\theta_{j-1}: h_i^*\neq h^*_{\bar\theta_{j-1}}\}$
to denote the successive jump time of the walk $h_i^*$,
one can write
$$\eta_i^*=\sum_{j:\bar \theta_j<\eta_i^*} G_j$$
where the $G_j$ are independent geometric random variables
with parameter $(\lambda+1)/(\lambda+
C_1\log t_{\epsilon})$ that represent the holding times.
Therefore, for any constants $C_2,C_3$ independent of $\epsilon$ and $t$,
$$P(\eta^*_{t_{\epsilon}}>C_2 t_{\epsilon} (\log t_{\epsilon})^2)
\leq P(\bar \theta_{C_3 t_{\epsilon}}< \eta^*_{t_{\epsilon}}
)+P(
\sum_{j=1}^{C_3 t_{\epsilon}} G_j>C_2 t_{\epsilon}(\log t_{\epsilon})^2)\,.$$
The event
$\{\bar\theta_{C_3 t_{\epsilon}}<\eta^*_{t_{\epsilon}}\}$
has the same probability as 
the event that a biased nearest neighbor random walk on $\mathbb{Z}$ started
at $0$,
with probability $\lambda/(\lambda+1)$ to increase at each step, 
does not hit
$t_{\epsilon}$ by time $C_3 t_{\epsilon}$. Because $\lambda>1$,
choosing $C_3=C_3(\epsilon)$ 
large,
this probability can be made exponentially small in $t_{\epsilon}$, and in 
particular bounded above by $1/t$ for $t$ large.
Fix such a $C_3$. Now,
$$P(
\sum_{j=1}^{C_3 t_{\epsilon}} G_j>C_2 t_{\epsilon}(\log t_{\epsilon})^2)
\leq C_3 t_{\epsilon} P(G_1>C_2 (\log t_{\epsilon})^2/C_3)\,.$$
By choosing $C_2=C_2(\epsilon)$ large, one can make this last term smaller
than $1/t$. Therefore,
with such a choice of $C_2$ and $C_3$,
and writing 
$\widehat{ \Upsilon}_{t, \epsilon}=\Upsilon_{t,\epsilon}\cap 
\{\eta_{t_{\epsilon}}<C_2 t_{\epsilon} (\log t_{\epsilon})^2\}$,
we obtain
	from (\ref{eq-150206ec}) that for all $t$ large,
\begin{equation}
	\label{eq-150206ecc}
	E^o_{\IGW}(N_t(\alpha))\leq  
4+
E^o_{\IGW}(N_t(\alpha); \xi_{t_{\epsilon}}\geq t; \widehat{
\Upsilon}_{t,\epsilon})\,. 
\end{equation}
On the event 
$\Upsilon_{t,\epsilon}$, all excursions 
$\{X_\ell,\; \ell=\eta_{i-1},\ldots,\eta_i-1\}$
away from $\Ray$
that start at $v\in \Ray$ with $h(v)>-t_\epsilon$ are excursions into 
${GW}$-trees where the degree of the root is bounded by
$C_1(\log t_\epsilon)-1$. Therefore,
\begin{align}
	\label{eq-210306a}
&E^o_{\IGW}\left(
\sum_{\ell=\eta_{i-1}}^{\eta_{i}} 
{\bf 1}_{X_\ell\in  Q_t(\mathcal{T})}; 
\Upsilon_{t,\epsilon}, 
h(X_{\eta_{i-1}})>-t_\epsilon
\right) \\
& \leq
\max_{d\leq C_1 (\log t_{\epsilon})-1}
E^o_{\widehat{\GW}}\left( 
\sum_{\ell=0}^{T_o}{\bf 1}_{h(X_\ell)\leq t^\alpha}| d_o=d\right)
\nonumber \\
& = 
\max_{d\leq C_1 (\log t_{\epsilon})-1}
\left(\sum_{j=0}^{t^\alpha}
E^o_{\widehat{\GW}}( 
{\cal N}_o(j)|d_o=d)\right)\,.
\nonumber \end{align}
Therefore,  for all $t$ large,
\begin{eqnarray}
&&E^o_{\IGW}(N_t(\alpha); \xi_{t_{\epsilon}}\geq t; \widehat{
\Upsilon}_{t,\epsilon}) 
)\nonumber\\
&\leq&
E^o_{\IGW}\left(
\sum_{i=1}^{C_2 t_{\epsilon}(\log t_{\epsilon})^2}
{\bf 1}_{\{h(X_{\eta_{i-1}})>-t_\epsilon\}}
\sum_{\ell=\eta_{i-1}}^{\eta_{i}} 
{\bf 1}_{X_\ell\in  Q_t(\mathcal{T})}; \Upsilon_{t,\epsilon} \right)
\nonumber\\
&\leq &  {C_2 t_{\epsilon}(\log t_{\epsilon})^2}  
\max_{d\leq C_1 (\log t_{\epsilon})-1}
\left(\sum_{j=0}^{t^\alpha}
E^o_{\widehat{\GW}}( 
{\cal N}_o(j)|d_o=d)\right)
\nonumber\\
&\leq &{t^{1/2+\alpha+\epsilon/2}}\,,
\end{eqnarray}
where the second inequality uses (\ref{eq-210306a}),
and
(\ref{eq-071005c})
was used in the last inequality.
Combined with (\ref{eq-150206ecc}), this completes the proof of 
Lemma \ref{lem-071005aa}.
\qed

\begin{corollary}
\label{cor-101005a}
For each $\epsilon>0$
there exists a $t_1=t_1({\cal T},\epsilon)<
\infty$ such that for all $t\geq t_1$,
\begin{equation}
\label{eq-101005a}
E_{\cal T}^oN_t(\alpha)\leq  
t^{1/2+\alpha+2\epsilon}\,, \; \IGW-a.s.\,.
\end{equation}
\end{corollary}
{\bf Proof of Corollary \ref{cor-101005a}}
From Lemma \ref{lem-071005aa} and Markov's inequality
we have
$$P_{\IGW}(E_{\cal T}^o N_t(\alpha)>
c_{\epsilon} t^{1/2+\alpha+3\epsilon/2})\leq t^{-\epsilon/2}\,.$$
Therefore, with $t_k=2^k$, it follows from  Borel-Cantelli
that there exists an $k_1=k_1({\cal T},\epsilon)$ such that
for $k>k_1$,
$$E_{\cal T}^oN_{t_k}(\alpha)\leq c_{\epsilon} 
t_k^{1/2+\alpha+3\epsilon/2}
\,,\; \IGW-a.s.\,.$$
But for $t_k<t<t_{k+1}$ one has that $N_t(\alpha)\leq N_{t_{k+1}}(\alpha)$.
The claim follows.
\qed

\noindent
{\bf Proof of Proposition \ref{prop-101005a}}
Note that the number of visits of $X_i$ to $Q_t(\mathcal{T})$
between time $i=t$ and $i=t+\lceil t^\delta\rceil$ is bounded
by $N_{t+\lceil t^\delta \rceil}(\alpha)$. Therefore,
$$P_{\cal T}^o(
X_{\tau_t}\in Q_t(\mathcal{T}))=
\frac{1}{t^\delta}\sum_{i=t}^{t+\lceil t^\delta\rceil} 
P_{\cal T}^o\left(X_i
\in Q_t(\mathcal{T})\right)
\leq
\frac{1}{t^\delta} E_{\cal T}^o(N_{t+\lceil t^\delta\rceil}(\alpha))\,.$$
Applying Corollary 
\ref{cor-101005a}  with our choice of
$\epsilon_0$, see 
(\ref{eq-150206}), it follows that
for all $t>t_1({\cal T},\epsilon_0)$, for $\IGW$-almost every $\cal T$,
$$P_{\cal T}^o(
X_{\tau_t}\in Q_t(\mathcal{T}))
\leq  \frac{(t+\lceil t^\delta \rceil)^{1/2+\alpha+3\epsilon_0}}
{t^{\delta}}\leq \frac{1}{t^{\epsilon_0}}\,.$$
\qed

\section{From \IGW \ to \GW: Proof of Theorem \ref{theo-GW}}
\label{sec-GWproof}
Our proof of Theorem \ref{theo-GW} is based on constructing 
a shifted coupling between the random walk $\{X_n\}$ on 
a \GW \  tree and a random walk $\{Y_n\}$ on an $\IGW$ tree. We begin by
introducing notation. For a tree (finite or infinite, rooted
or not) ${\cal T}$, we let 
$ {\cal LT}$ denote the collection of leaves
of ${\cal T}$, that is of  vertices of degree $1$ in ${\cal T}$
other than the root.
We set ${\cal T}^o={\cal T}\setminus {\cal LT}$.
For two trees ${\cal T}_1,{\cal T}_2$ with roots (finite or infinite) 
and a vertex $v\in {\cal LT}_1$, we let
${\cal T}_1\circ ^v {\cal T}_2$ denote the tree obtained 
by gluing the root of ${\cal T}_2$ at the vertex $v$.
Note that if ${\cal T}_1$ has an infinite ray emanating
from the root, and ${\cal T}_2$ is a finite rooted tree, then
${\cal T}_1\circ ^v {\cal T}_2$ is a rooted tree with a marked 
infinite ray emanating from the root.

Given a \GW\ tree ${\cal T}$ and a path $\{X_n\}$ on the tree, 
we construct a family of finite trees ${\cal T}_i$ and 
of finite paths $\{u_n^i\}$ on ${\cal T}_i$ as follows.
Set $\tau_0=0,\eta_0=0$, and let ${\cal U}_0$ denote the rooted
tree consisting of the root $o$ and its offspring.
For $i\geq 1$, let
\begin{eqnarray}
\label{eq-051105a}
\tau_i&=&\min\{n>\eta_{i-1}: X_n\in  {\cal LU}_{i-1}\},
\quad (\mbox{\rm Excursion start})\nonumber \\
\eta_i&=&\min\{n>\tau_i: X_n \in  {\cal U}_{i-1}^o\},
\quad (\mbox{\rm Excursion end})\nonumber \\
v_i&=& X_{\tau_i}, \quad (\mbox{\rm Excursion start location})\,.
\end{eqnarray}
We then set
$${\cal V}_i= \{v\in {\cal T}: X_n=v \, \mbox{\rm for some
$n\in [\tau_i,\eta_i)$}\}\,,
$$
define $\overline{\cal V}_i={\cal V}_i \cup
\{v\in {\cal T}: \mbox{\rm $v$ is an offspring of some
$w\in {\cal V}_i$}\}$ and let ${\cal T}_i$ denote the rooted
subtree
of ${\cal T}$ with vertices in $\overline{\cal V}_i$ and
root $v_i$. We also define the path 
$\{u_n^i\}_{n=0}^{\eta_i-\tau_i-1}$ by
$u_n^i=X_{n+\tau_i}$, noting that $u_n^i$ is
a path in ${\cal T}_i$.
Finally, we set
\begin{equation}
\label{eq-051105b}
{\cal U}_i={\cal U}_{i-1}\circ ^{v_i}{\cal T}_i\,.
\end{equation}
Note that ${\cal U}_i$ is a  tree rooted at $o$ since $v_i\in 
{\cal LU}_{i-1}$. Further, by the $\GW$-almost sure
recurrence of the biased random walk
on ${\cal T}$, it holds that ${\cal T}=\lim_i {\cal U}_i$.

Next, we construct an $\IGW$ tree $\widehat {\cal T}$ with
root $o$ and an infinite ray, denoted $\Ray$, emanating
from the root, and a 
($\lambda$-biased) random walk $\{Y_n\}$ on $\widehat {\cal T}$,
as follows. First, we choose a vertex denoted $o$ and 
a semi-infinite directed path
$\Ray$ emanating from it. Next, we let each vertex $v\in
\Ray$ have 
$d_v$ offspring, where $P(d_v=k)=kp_k/m$, and the $\{d_v\}_{v\in
\Ray}$ are independent. For each vertex $v\in \Ray$, $v\neq o$,
we identify one
of its offspring with the vertex $w\in \Ray$ that satisfies
$d(w,o)=d(v,o)-1$, and write $\widehat {\cal U}_0$ for the resulting 
tree with root $o$ and marked ray $\Ray$. 

Set next $\hat \tau_0=
\hat \eta_0=0$.
We start a $\lambda$-biased random walk $Y_n$ on
$\widehat {\cal U}_0$ with $Y_0=o$, and define
$$ \hat \tau_1=\min \{n>0: Y_n\in {\cal L}\widehat {\cal U}_0\}\,.$$ 
Let $\hat v_1=Y_{\hat \tau_1}$. We now set
$\widehat {\cal U}_1=\widehat {\cal U}_0\circ^{\hat v_1} {\cal T}_1$ and
$\hat \eta_1=\hat\tau_1+\eta_1-\tau_1$, and for
$\hat \tau_1\leq n\leq \hat \eta_1-1$,
set $Y_n=u_{n-\hat\tau_1}^i$. Finally, with $\hat w_1$ the ancestor of
$\hat v_1$, we set $Y_{\hat\eta_1}=\hat w_1$.

The rest of the construction proceeds similarly. For $i>1$,
start a $\lambda$-biased 
random walk $\{Y_n\}_{n\geq \hat\eta_{i-1}} $ 
on $\widehat {\cal U}_{i-1}$ with $Y_{\hat \eta_{i-1}}= \hat 
w_{i-1}$ and define 
\begin{eqnarray}
\label{eq-051105c}
\hat \tau_i&=&\min\{n>\hat \eta_{i-1}: Y_n\in {\cal L}\widehat 
{\cal U}_{i-1}\},
\quad (\mbox{\rm Excursion start})\,,\nonumber \\
\hat v_i&=& Y_{\hat
\tau_i}, \quad (\mbox{\rm Excursion start location})\,,\\
\hat \eta_i&=&\hat\tau_i+\eta_i-\tau_i\,,
\quad (\mbox{\rm Excursion end})\,,\nonumber \\
\widehat {\cal U}_i&=& \widehat {\cal U}_{i-1}\circ^{\hat v_i}{\cal T}_i\,,
\quad (\mbox{\rm Extended  tree}),\nonumber\\
Y_n&=& X_{n-\hat\tau_i}, n\in [\hat\tau_i,\hat \eta_i)
\quad (\mbox{\rm Random walk path during excursion})\,,\nonumber\\
Y_{\hat \eta_i}&=&\hat w_i= 
\mbox{\rm ancestor of}\; \hat v_i\,.\nonumber
\end{eqnarray}
Finally, with $\widehat {\cal U}=\lim_i \widehat {\cal U}_i$, define
the tree $\widehat {\cal T}$ by attaching to each vertex of ${\cal L}
\widehat {\cal U}$ an independent Galton-Watson tree, thus obtaining an 
infinite tree with root $o$ and infinite ray
emanating from it. 
The construction leads immediately to
the following.
\begin{lemma}
\label{lem-coupling}
a) The tree $\widehat {\cal T}$ with root $o$ and marked ray $\Ray$
is distributed according to $\IGW$.\\
b) Conditioned on $\widehat {\cal T}$,
the law of 
$\{Y_n\}$ is the law of a $\lambda$-biased 
random walk on $\widehat{\cal T}$.\\
\end{lemma}
\begin{figure}[htp]
\psfrag{t=1}{$t=1$}
\psfrag{(tau1=1)}{$(\tau_1=1)$}
\psfrag{t=6}{$t=6$}
\psfrag{(eta1=6)}{$(\eta_1=6)$}
\psfrag{calU1}{${\cal U}_1$}
\psfrag{t=2}{$t=2$}
\psfrag{(t=7)}{$t=7$}
\psfrag{hatcalU1}{$\hat{\cal U}_1$}
\psfrag{0}{$\mathbf{o}$}
\subfiguretopcaptrue
    \subfigure[GW side.]{\label{fig:GWside}\epsfig{file=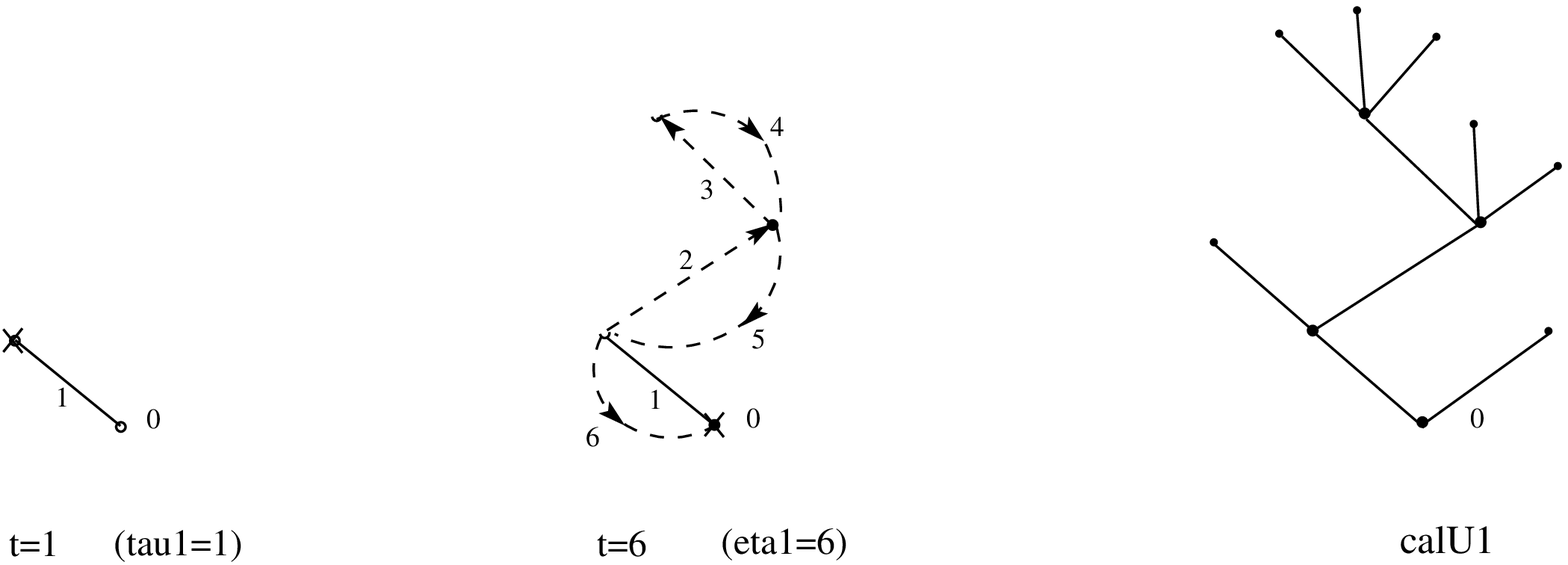,height=2.1in, width=6in, angle=0}} \\
\vspace{0.5in}%
\subfiguretopcaptrue
    \subfigure[IGW side.]{\label{fig:IGWside}\epsfig{file=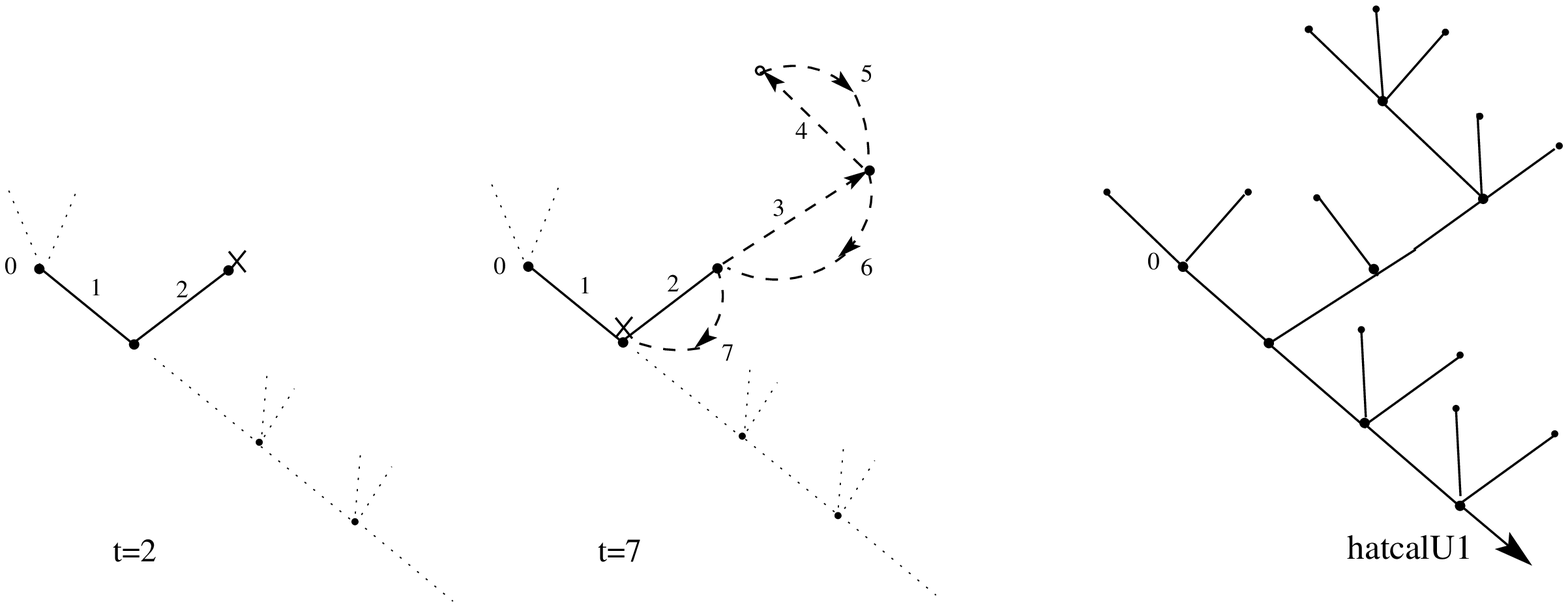,height=2.1in, width=6in, angle=0}} 
 \begin{center}
\caption{The coupling between the \GW \ and \IGW \ walks. 
X marks the location of the walker.}
 \end{center}
\end{figure}
Let ${\cal R}_n=h(Y_n)-\min_{i=1}^n h(Y_i)\geq 0$.
Due to Theorem
\ref{theo-1}, for $\IGW$-almost all $\widehat {\cal T}$,
the process ${\cal R}_{\lfloor nt\rfloor}/\sqrt{n}$
converges to a Brownian motion reflected at its running minimum,
which possesses the same law as  the absolute value 
of a Brownian motion,  see e.g. \cite[Theorem 6.17]{KS}. Our efforts
are therefore directed toward estimating the relation between
the processes $\{X_n\}$ and $\{{\cal R}_n\}$.
Toward this end, 
let $I_n=\max \{i: \tau_i\leq n\}$ and $\widehat I_n=\max\{i: \hat\tau_i
\leq n\}$ measure the number of excursions started by the
walks $\{X_n\}$ and $\{Y_n\}$ before time $n$, and set
$\Delta_n=\sum_{i=1}^{I_n} (\tau_i-\eta_{i-1})$,
and $\widehat 
\Delta_n=\sum_{i=1}^{\hat I_n} (\hat \tau_i-\hat \eta_{i-1})$.
Set also
$B_n=\max_{s<t\leq n: Y_s\in \Ray, Y_t\in \Ray}(h(Y_t)-h(Y_s))$
($B_n$ measures the maximal amount the random walk $\{Y_n\}$
backtracks, that is moves against the drift,
along $\Ray$ before time $n$).
Next set, 
 recalling (\ref{eq-190106b}),
\begin{eqnarray}
\label{eq-051105h}
\Delta_n^\alpha&=& \sum_{i=1}^{I_n} 
\sum_{t\in [\eta_{i-1},\tau_i)} {\bf 1}_{|X_t|\leq n^{\alpha}}\,,
\nonumber \\
\widehat\Delta_n^\alpha&=& \sum_{i=1}^{\widehat I_n} 
\sum_{t\in [\hat \eta_{i-1},\hat \tau_i)} {\bf 1}_{Y_t\in
Q_{n^\alpha}(\widehat {\cal T})}\,.
\end{eqnarray}
Clearly, $\Delta_n^\alpha\leq \Delta_n$
and $\widehat \Delta_n^\alpha\leq \widehat \Delta_n$. We however
can say more.
\begin{lemma}
\label{lem-dn}
Let $A_n=\{\Delta_n^\alpha= \Delta_n\}$ and $\widehat A_n=\{\widehat 
\Delta_n^\alpha=
\widehat \Delta_n\}$. Then,
\begin{equation}
\label{eq-051105g1}
\lim_{n\to\infty}
P_{\cal T}^o(A_n^c)=0\,,\quad \GW-a.s.,
\end{equation}
and
\begin{equation}
\label{eq-051105g2}
\lim_{n\to\infty}
P_{\cal T}^o(\widehat{A_n}^c)=0\,,\quad \IGW-a.s..
\end{equation}
Further,
\begin{equation}
\label{eq-051105g3}
\limsup \frac{|\Delta_n|}{n}=0\,,\quad \GW-a.s.\,,
\end{equation}
and
\begin{equation}
\label{eq-051105g4}
\limsup \frac{|\widehat \Delta_n|}{n}=0\,,\quad \IGW-a.s.\,,
\end{equation}
Finally,
\begin{equation}
\label{eq-051105g5}
\limsup \frac{B_n}{\sqrt n}=0\,,\quad \IGW-a.s.\,,
\end{equation}
\end{lemma}
We postpone for the moment the proof of Lemma \ref{lem-dn}.
Note that on  the event $A_n\cap \widehat A_n$,
one has 
\begin{equation}
\label{eq-051105e}
\min_{s: |s-n|\leq \Delta _n+ \widehat \Delta_n}
|\; |X_n|-{\cal R}_{s}|\leq  2n^\alpha+B_n\,.
\end{equation}
(To see that, note that the position $|X_n|$ consists of sums of excursions
$\{u_\cdot^i\}$, up to an error coming from the
parts of the path not contained in these excursions, all
contained in a distance at most $n^\alpha$ from the root. Similarly, for
some $s$ with $|s-n|\leq \Delta_n+\widetilde{\Delta}_n$,
${\cal R}_s$ consists of the sum of the same excursions, up to an error
coming from the parts of the path not contained in these excursions, which 
sum up to a total distance of at most $n^\alpha$ from $\Ray$ in addition
to the amount $B_n$ of backtracking along $\Ray$.)

In view of Lemma \ref{lem-dn}, the convergence in distribution
(for $\IGW$-almost every $\widehat {\cal T}$) of
${\cal R}_{\lfloor n t\rfloor}/\sqrt{n}$ to reflected Brownian
motion, together with 
(\ref{eq-051105e}),
complete the proof of Theorem \ref{theo-GW}.
\qed

\noindent
{\bf Proof of Lemma \ref{lem-dn}}
Consider a rooted tree $\cal T$ distributed according
to $\GW$, and a random walk path $\{X_t\}_{t\geq 0}$ with $X_0=o$
on it.
We introduce some  notation. For $k\geq 1$, 
let $a_k=\sum_{j=1}^k \tau_j$, 
$b_k=\sum_{j=1}^{k-1} \eta_j$, and $J_k=[a_k-b_k+k,a_{k+1}-b_{k+1}+k]$
(the length of $J_k$ is the time spent by the walk between the 
$k$-th and the $k+1$-th excursions). For $s\in J_k$, we define 
$t(s)=\eta_k+s-(a_k-b_k+k)$. Finally, we set $\widetilde{X}_0=0,
\widetilde{X}_1=
X_{\tau_1}=X_1$, and $\widetilde{X}_s=X_{t(s)}$ (note that
the process $\widetilde{X}_s$ travels on  vertices
``off the coupled excursions''). Note that even conditioned on
${\cal T}$, the nearest neighbor process 
$\{\widetilde{X}_s\}_{s\geq 0}$ on ${\cal T}$ is neither Markovian nor 
progressively measurable
with respect to 
its natural filtration. To somewhat address this issue, we define the
filtration
${\cal G}_s=\sigma(X_i, i\leq t(s))\,,$ and note
that conditioned on ${\cal T}$,
$\{\widetilde{X}_s\}_{s\geq 0}$ is progressively measurable with respect to the
filtration ${\cal G}_s$.

The statement (\ref{eq-051105g1})
will follow as soon as we prove the statement
\begin{equation}
	\label{eq-170506a}
\lim_{n\to\infty}
P_{\cal T}^o(\max_{s\in\cup_{k=1}^{I_n} J_k}
|\widetilde{X}_s|\geq n^\alpha)=0\,,
\quad \GW-a.s.,
\end{equation}
The proof of (\ref{eq-170506a}) will be carried out in several steps. The
first step allows us to control the event that the
time spent by the process $X_t$ 
inside excursions is short. The proof is routine and postponed.
\begin{lemma}
	\label{lem-170506b}
	For all $\epsilon>0$, 
\begin{equation}
		\label{eq-170605c}
\lim_{n\to\infty}
		P_{\cal T}^o(\sum_{i=1}^{n^{1/2+\epsilon}}
		(\eta_i-\tau_i)<n)=0\,,\quad \GW-a.s..
	\end{equation}
Further, with 
$$\widetilde{T}_n=\min\{t: W_{X_t}>(\log n)^2\}\,,$$
it holds that
\begin{equation}
	\label{eq-170506f}
	\lim_{n\to\infty}
	n P_{\cal T}^o(\widetilde{T}_n\leq n)=0\,,\quad \GW-a.s.
\end{equation}
\end{lemma}
Our next step involves ``coarsening'' the process $\{\widetilde{X}_s\}$
by stopping it at  random times $\{\Theta_i\}$ in such a way that if the 
stopped process has increased its distance from the root between two 
consecutive stopping times, with high probability one of the intervals
$J_k$ has been covered. More precisely, define
$\Theta_0=0$, and for $i\geq 1$, 
$$\Theta_i=\min\{s>\Theta_{i-1}: \left|
|\widetilde{X}_s|-|\widetilde{X}_{\Theta_{i-1}}|\right|=
\lfloor (\log n)^{3/2}\rfloor\}\,.$$
We emphasize that the $\Theta_i$ depend on $n$, although this dependence
is suppressed in the notation.
The following lemma, whose proof is 
again routine and postponed, explains why this coarsening is useful.
\begin{lemma}
	\label{lem-170506h}
	With the notation above, 
	\begin{equation}
		\label{eq-170506j}
		\lim_{n\to\infty}
		P_{\cal T}^o(\mbox{\rm for some $k\leq I_n$,
		$\Theta_{i-1},\Theta_i\in J_k$},\,
		|\widetilde{X}_{\Theta_i}|>|\widetilde{X}_{\Theta_{i-1}}|
		)=0\,,\quad 
		\GW-a.s.
	\end{equation}
\end{lemma}
We have now prepared all needed preliminary steps. Fix $\epsilon>0$.
Note first that
due to (\ref{eq-oof}) and the Borel-Cantelli lemma,
for all $n$ large,
$|A_{n^\alpha}^\epsilon|\leq |D_{n^\alpha}|e^{-\nu(\epsilon) n^\alpha}$,
\GW-a.s. On the other hand, since $E_{\GW} |D_{n^\alpha}|=m^{n^\alpha}$,
Markov's inequality and the Borel-Cantelli lemma imply that for all $n$ large,
$|D_{n^\alpha}|\leq m^{n^{\alpha}}e^{\nu(\epsilon) n^\alpha/2}$,
\GW-a.s. Combining these facts, it holds that for all $n$ large, 
\begin{equation}
\label{eq-16-6-6b}
|A_{n^\alpha}^\epsilon|\leq m^{n^\alpha} e^{-\nu(\epsilon) n^\alpha/2},
\quad
\GW-a.s.\,.
\end{equation} 
For any vertex $v\in D_{n^\alpha}$, by considering the trace of the 
random walk on the path connecting $o$ and $v$ it follows that
$$P_{\cal T}^o (X_t=v \; \mbox{\rm for some $t\leq n$})
\leq  1-(1-\lambda^{-n^{\alpha}})^n\leq n\lambda^{-n^{\alpha}}\,,
\GW-a.s.$$
 Using this and (\ref{eq-16-6-6b}) in the first inequality, and 
(\ref{eq-170506f}) in the second, we get
\begin{eqnarray}
\label{eq-160606d}
&&\limsup_{n\to\infty}
P_{\cal T}^o(\max_{s\in\cup_{k=1}^{I_n} J_k}
|\widetilde{X}_s|\geq n^\alpha)\\
&\leq &
\limsup_{n\to\infty}
	P_{\cal T}^o(\exists s\in\cup_{k=1}^{I_n} J_k:
|\widetilde{X}_s|= n^\alpha, S_{\widetilde{X}_s}\geq \eta n^{\alpha}/2)
\nonumber \\
&\leq &
\limsup_{n\to\infty}
	P_{\cal T}^o(\exists s\in\cup_{k=1}^{I_n} J_k:
|\widetilde{X}_s|= n^\alpha, S_{\widetilde{X}_s}\geq \eta n^{\alpha}/2,
t(s)\leq \widetilde{T}_n)\,.\nonumber
\end{eqnarray}
We next note that by construction,
$$ |\{i\in \{1,\ldots,\ell\}: |\widetilde{X}_{\Theta_i}|>
|\widetilde{X}_{\Theta_{i-1}}|\}|\geq \ell/2\,.$$
Hence, with $P_{\cal T}$ probability approaching $1$ as $n$
goes to infinity, $t(\Theta_{2n^{1/2+\epsilon}})>n$ because of 
(\ref{eq-170605c})
and Lemma \ref{lem-170506h}.
From this and (\ref{eq-160606d}),
we conclude that
\begin{eqnarray*}
&&\limsup_{n\to\infty}
	P_{\cal T}^o(\max_{s\in\cup_{k=1}^{I_n} J_k}
|\widetilde{X}_s|\geq n^\alpha)\\
&\leq& 
	\limsup_{n\to\infty}
	\sum_{i=1}^{2n^{1/2+\epsilon}}
	P_{\cal T}^o\Big(
|\widetilde{X}_{\Theta_i}|\geq n^\alpha-(\log n)^2 , 
S_{\widetilde{X}_{\Theta_i}}\geq \eta n^{\alpha}/2-(\log n)^4,\\
&&\quad \quad \quad \quad
\quad \quad \quad \quad
\quad \quad \quad \quad
\widetilde{T}_n
>t(\Theta_{i})\Big)\,.
\end{eqnarray*}
On the event $\widetilde{T}_n
>t(\Theta_{i})$ it holds that 
$|S_{\widetilde{X}_{\Theta_{i}}}-S_{\widetilde{X}_{\Theta_{i-1}}}|
\leq (\log n)^4\;.$ Therefore,
decomposing according to return
times of $\widetilde{X}_{\Theta_i}$ to the root,
\begin{eqnarray}
	\label{eq-170506k}
	&&\limsup_{n\to\infty}
	P_{\cal T}^o(\max_{s\in\cup_{k=1}^{I_n} J_k}
|\widetilde{X}_s|\geq n^\alpha)\\
&\leq &
\limsup_{n\to\infty} 
\sum_{i=0}^{2n^{1/2+\epsilon}}
\sum_{j=i+1}^{2n^{1/2+\epsilon}}
P_{\cal T}^o\Big(
|\widetilde{X}_{\Theta_j}|\geq n^\alpha-(\log n)^2 , 
\widetilde{X}_{\Theta_i}=o, \nonumber \\
&&
\quad \quad \quad \quad 
S_{\widetilde{X}_{\Theta_j}}\geq \eta n^{\alpha}/2-(\log n)^4,\nonumber
\\
&&
\quad \quad \quad \quad 
|\widetilde{X}_{\Theta_k}|>0\, 
\mbox{\rm and}
\;
|S_{\widetilde{X}_{\Theta_{k}}}-S_{\widetilde{X}_{\Theta_{k-1}}}|
\leq (\log n)^4\;
\mbox{\rm for $i < k\leq j$}\,
\Big)\,\nonumber\\
&=:& 
\limsup_{n\to\infty} 
\sum_{i=0}^{2n^{1/2+\epsilon}}
\sum_{j=i+1}^{2n^{1/2+\epsilon}}
P_{i,j,n}\,.
\nonumber\end{eqnarray}
Fixing $i$,
set for $t\geq 1$, 
$\widetilde{M}_t=
S_{\widetilde{X}_{\Theta_{i+t}}}$.
Introduce the random time 
\begin{eqnarray*}
&&K_n=\min\{t>1: 
\mbox{\rm  $X_s=o$ 
for some $s\in [t(\Theta_{i+1}),t(\Theta_{i+t})]$} \,\\
&&\quad \quad 
\quad \quad \quad \quad 
\quad \quad \quad \quad 
\mbox{\rm or
$|\widetilde{M}_{t}-\widetilde{M}_{t-1}|\geq (\log n)^4$}\}\,,
\end{eqnarray*}
and the filtration $\widetilde{\cal G}_t={\cal G}_{\Theta_{i+t}}$.
The crucial observation is that 
$\{\widetilde{M}_{t\wedge K_n}-\widetilde{M}_1\}$
is a supermartingale for the filtration $\widetilde{\cal G}_t$, 
with increments bounded in absolute value by $(\log n)^4$ for all $t<K_n$,
and bounded below by $-(\log n)^4$ even for $t=K_n$
(it fails to be a martingale due to the ``defects'' 
at the boundary of each of the intervals $J_k$, at which times $r$ the 
conditional expectation of the increment $S_{\widetilde{X}_{r+1}}-
S_{\widetilde{X}_r}$
is negative). 
Let $\widetilde{M}_t'=\widetilde{M}_t$ if $t<K_n$ or $t=K_n$
but $\widetilde{M}_{t}<\widetilde{M}_{t-1}+(\log n)^4$,
and $\widetilde{M}_t'=\widetilde{M}_{K_n-1}$ otherwise.
That is, $\widetilde{M}'_t-\widetilde{M}_1$ 
is a truncated version of the 
supermartingale $\widetilde{M}_{t\wedge K_n}
-\widetilde{M}_1$. It follows that
for some non-negative process $a_t$,
$\{\widetilde{M}_{t}'-\widetilde{M}_1+a_t\}$
is a martingale with increments bounded for all $t\leq K_n$ by $2(\log n)^4$.
Therefore, by Azuma's inequality \cite{azuma},
for $j\leq n^{1/2+\epsilon}$, and all $n$ large,
$$P_{i,j,n}\leq
P_{\cal T}^o\Big(\max_{1\leq k\leq 2n^{1/2+\epsilon}}
[\widetilde{M}_{k}'-\widetilde{M}_1]
\geq \eta n^\alpha/3\Big)
\leq e^{-n^{2\alpha}/n^{1+3\epsilon}}\,.
$$
Since this estimate did not depend on $i$ or $j$,
together with (\ref{eq-170506k}), 
this completes the proof of (\ref{eq-170506a}), and hence of
(\ref{eq-051105g1}). The proofs of (\ref{eq-051105g2}) 
and
(\ref{eq-051105g5})
are similar and 
omitted.


We next turn to the proof of (\ref{eq-051105g4}). 
Recall that from Lemma \ref{lem-071005aa}, for 
any $\epsilon>0$, and all $n>n_0(\epsilon)$, 
$$P_{\IGW}^o(N_n(\alpha)\geq  n^{1/2+\alpha+2\epsilon})
\leq n^{-\epsilon}\,.$$
Therefore, noting the monotonicity of
$N_n(\alpha)$ in $n$, an application of the Borel-Cantelli
lemma (to the sequence $n_k=2^k$) shows that
$$ \frac{N_n(\alpha)}{n^{1/2+\alpha+3\epsilon}}\to_{n\to \infty}
 0\,, \quad \IGW-a.s.$$
Since $\epsilon$ can be chosen such that
$1/2+\alpha+3\epsilon<1$, c.f. (\ref{eq-150206}),
and  $\widehat \Delta_n^\alpha\leq N_n(\alpha)$,
(\ref{eq-051105g4}) follows.

We finally turn to the proof of 
(\ref{eq-051105g3}). In what follows, we let
$C_i=C_i({\cal T})$ denote  constants that may depend
on ${\cal T}$ (but not on $n$).
Let
$T_{\epsilon}(n)=\min\{t: |X_t|=n^{1/2+\epsilon}\}$.
By Lemma \ref{lem-111005a},
$$P_{\GW}^o(T_{\epsilon}(n)\leq n)\leq 4 n e^{-n^{2\epsilon}/2}\,.$$
In particular, 
by the Borel-Cantelli lemma,  for $\GW$-almost every ${\cal T}$,
\begin{equation}
\label{eq-051105oof}
P_{\cal T}^o(T_{\epsilon}(n)\leq n)\leq C_4({\cal T})
 e^{-n^{\epsilon}}\,.
\end{equation}
Let ${\cal C}_{o,\ell}$ denote the  conductance between the root and $D_\ell$.
That is, define a unit flow $f$ on ${\cal T}$ as a collection
of non-negative numbers $f_{v,w}$, with $v\in \cal T$ and $w
\in {\cal T}$
an offspring of $v$, such that Kirchoff's current law hold:
$1=\sum_{w\in D_1} f_{o,w}$
and $f_{v,w}=\sum_{w': w'\, \mbox{\rm  is an  offspring of }\, w}
f_{w,w'}$. Then,
$$ {\cal C}_{o,\ell}^{-1}
=\inf_{f: f \,\mbox{\rm \small is a unit flow}}
\sum_{i=0}^{\ell-1} \sum_{v\in D_i}
\sum_{w: w \,\mbox{\rm \small is an offspring of}\, v}
f_{v,w}^2 \lambda^i\,.$$
By \cite[Theorem 2.2]{PP},  for $\GW$-almost every ${\cal T}$
there exists  a constant $C_5({\cal T})$ and a unit flow $f$
such that
$$\sum_{v\in D_i}
\sum_{w: w\, \mbox{\rm \small is an offspring of}\, v}
f_{v,w}^2\leq C_5({\cal T}) \lambda^{-i}\,.$$
It follows that 
\begin{equation}
	\label{eq-180506a}
	{\cal C}_{o,\ell}^{-1}\leq C_5({\cal T}) \ell.
\end{equation}
On the other hand,
by standard theory, see \cite[Exercise 2.47]{LP},
for a given tree
$\cal T$, with $L_o(j)$ denoting the number
of visits to the root before time $j$,
$$E^o_{\cal T} L_o(T_{\epsilon}(n))=
d_o
{\cal C}_{o,n^{1/2+\epsilon}}^{-1}\,.$$
Hence, 
$E^o_{\cal T} L_o(T_{\epsilon}(n))\leq d_oC_5({\cal T}) n^{1/2+\epsilon}$.
By Lemma \ref{lem-071005b}, 
we also have that\\
$E^o_{\cal T}({\cal N}_o(\ell))\leq C_6({\cal T})$, for any $\ell$.
Thus, using $\bar N_n(\alpha)=\sum_{t=1}^n {\bf 1}_{\{|X_t|\leq n^\alpha\}}$,
$$ E^o_{\cal T}( \bar N_n(\alpha);
T_{\epsilon}(n)\geq n)
\leq 
E^o_{\cal T} L_o(T_{\epsilon}(n))
E^o_{\cal T}\!\! \left(\sum_{\ell=0}^{n^\alpha} {\cal N}_o(\ell)\right)
\!\!\leq
d_o C_5({\cal T})C_6({\cal T}) n^{1/2+\epsilon+\alpha}\,.$$
It follows from this that
$$ E^o_{\cal T}(\bar N_n(\alpha))\leq 
n P_{\cal T}^o(T_{\epsilon}(n)\leq n)+
d_oC_5({\cal T})C_6({\cal T}) n^{1/2+\epsilon+\alpha}\,.$$
Using (\ref{eq-051105oof}) and the fact that
$\bar N_n(\alpha)\geq \Delta_n^\alpha$, together with (\ref{eq-051105g1}),
completes the proof
of 
(\ref{eq-051105g3}), and hence of Lemma \ref{lem-dn}.
\qed

\noindent
{\bf Proof of Lemma \ref{lem-170506b}:} We note first that under the 
annealed measure $\GW$, the random times $(\eta_i-\tau_i)$, which denote
the length of the  excursions, are i.i.d., and
for all $x$,
$$P^o_{\GW}(\eta_i-\tau_i\geq x)\geq \frac{1}{\lambda+1}
P^o_{\GW}(T_o\geq x),$$
where $T_o=\min\{t\geq 1:X_t=o\}$ denotes the first return time of
$X_t$ to $o$.

Throughout, the constants 
$C_i({\cal T})$, that depend only 
on the tree ${\cal T}$, 
are as in the proof above.
Let $x_t= t^{1/2+\epsilon/2}$ and set $T_z=\min\{t: |X_t|=z\}$. Then,
\begin{equation}
	\label{tapasm1}
	P_{\cal T}^o(T_o\geq t)\geq P_{\cal T}^o(T_{x_t}<T_o)
P_{\cal T}^o(T_{x_t}\geq t|T_{x_t}<T_o)\,.
\end{equation}
Note however that 
$P_{\cal T}^o(T_{x_t}<T_o)$ is bounded by
the effective conductance between the root and $D_{x_t}$, which
by (\ref{eq-180506a}) is bounded below by 
$C_5({\cal T})
x_t^{-1} $.
In particular, 
\begin{equation}
	\label{tapas}
	P_{\cal T}^o(T_{x_t}<T_o)\geq \frac{C_5({\cal T})}{x_t}
\end{equation}
On the other hand, using (\ref{tapas}) and 
the Carne-Varopoulos bound 
(see \cite[Theorem 12.1]{LP},
\cite{carne,var}) in the second
inequality,
\begin{equation}
	\label{tapas-1}
	P_{\cal T}^o(T_{x_t}<t|T_{x_t}<T_o)
	\leq  \frac{P_{\cal T}^o(T_{x_t}<t)}{
P_{\cal T}^o( T_{x_t}<T_o)}
\leq  C_7({\cal T}) x_t e^{-t^{2\epsilon}}
\end{equation}
It follows that for all $t$ large,
$$P_{\cal T}^o(T_{x_t}\geq t|T_{x_t}<T_o)>1/2\,,$$
implying with (\ref{tapasm1}) and (\ref{tapas})  that for all $t$ large, 
\begin{equation}
	\label{tapas3}
	P_{\cal T}^o(T_o\geq t)\geq 
	\frac{C_5({\cal T})}
{2t^{1/2+\epsilon/2}}\,.
\end{equation}
It follows that for some deterministic constant $C$ and all $t$ large,
\begin{equation}
	\label{tapas4}
	P_{\GW}^o(T_o\geq t)\geq 
	\frac{C}
{t^{1/2+\epsilon/2}}\,.
\end{equation}
Hence,
\begin{eqnarray*}
	P_{\GW}^o(\sum_{i=1}^{n^{1/2+\epsilon}}
(\eta_i-\tau_i)<n) &\leq &
\left(1-\frac{P_{\GW}^o(T_o\geq n)}{\lambda+1}\right)^{n^{1/2+\epsilon}}
\\&\leq& \left(1-\frac{C}{n^{1/2+\epsilon/2}}\right)^{n^{1/2+\epsilon}}
\leq e^{-C n^{\epsilon/2}}\,.
\end{eqnarray*}
An application of the Borel-Cantelli lemma yields (\ref{eq-170605c}).

To see (\ref{eq-170506f}), note that by time $n$ the walker explored at most
$n$ distinct sites.
We say that $t$ is a fresh time if $X_s\neq X_t$ for all $s<t$. Then,
$$P_{\GW}^o(W_{X_t}\geq (\log n)^2, \, \mbox{\rm $t$ \, is a fresh time})
\leq P_{\GW}^o(W_o\geq (\log n)^2)\leq e^{-c(\log n)^2}\,,$$
by the tail estimates on $W_o$, see \cite{At94}. Therefore,
\begin{eqnarray*}
	&&P_{\GW}^o(W_{X_t}\geq (\log n)^2, \, \mbox{\rm for some
$t\leq n$})\\
&&\leq
\sum_{t=0}^n
P_{\GW}^o(W_{X_t}\geq (\log n)^2, \, \mbox{\rm $t$ \, is a fresh time})
\leq (n+1) e^{-c(\log n)^2}\,,
\end{eqnarray*}
from which (\ref{eq-170506f}) follows by an application of the Borel-Cantelli
lemma.\qed

\noindent
{\bf Proof of Lemma \ref{lem-170506h}:}
Let $G_n$ denote 
the event inside the probability in the left hand side of (\ref{eq-170506j}).
The event $G_n$
implies the existence of times $t_0<t_1<t_2\leq n$ and vertices
$u,v$ such that $X_{t_0}=u=X_{t_2}$, $X_{t_1}=v$, and
$|v|=|u|-\lfloor(\log n)^{3/2}\rfloor$. Thus, using the Markov property,
$$P^o_{\cal T}(G_n)\leq |\{(t_0,t_1):\,
t_0<t_1\leq n\}|
\max_{\stackrel{u,v\in {\cal T}:}{|v|=|u|-\lfloor(\log n)^{3/2}\rfloor
}}P^v(X_t=u \,
\mbox{\rm for some $t\leq n$})\,.$$
Noting that for each fixed $u,v$ as above, the last probability is
dominated by the probability of a $\lambda$-biased (toward $0$) 
random walk on $\mathbb{Z}_+$ reflected at $0$
to hit
location $\lfloor(\log n)^{3/2}\rfloor$ before time $n$, we get
$$P^o_{\cal T}(G_n)\leq n^2 e^{-c (\log n)^{3/2}}\,,$$
for some $c>0$, which implies (\ref{eq-170506j}).
\qed

\section{The transient case}
\label{sec-remarks}
Recall that when $\lambda<m$, it holds that $|X_n|/n\to_{n\to\infty} \v>0$,
$\GW$-a.s., for some non-random $\v=\v(\lambda)$ (see \cite{LPP2}).
Our goal in this section is to prove the following:
\begin{theorem}
\label{theo-GWtransient}
Assume $\lambda<m$ and $p_0=0, \sum_k \beta^k p_k<\infty$ for some
$\beta>1$. Then,
there exists a deterministic constant $\sigma^2>0$ such that for
$\GW$-almost every ${\cal T}$, the processes 
$\{(|X_{\lfloor nt \rfloor}|-nt\v)/\sqrt{\sigma^2 n}\}_{t\geq 0}$
converges in law to standard Brownian motion.
\end{theorem}

Before bringing the 
proof of Theorem \ref{theo-GWtransient}, we need to derived an {\it annealed}
invariance principle, see Corollary \ref{cor-190506g} below.
The proof of the latter proceeds
via the study of regeneration times, which are defined 
as follows: we set 
$$\tau_1:=\inf\{t: |X_t|>|X_s| \,\mbox{\rm for all $s<t$},
\,\mbox{\rm and}\, |X_u|\geq |X_t| \; \mbox{\rm for all}\; u\geq t\}\,,$$
and, for $i\geq 1$, 
$$\tau_{i+1}:=\inf\{t>\tau_i:
|X_t|>|X_s|\, \mbox{\rm for all $s<t$},\,
\mbox{\rm and}\, |X_u|\geq |X_t| \; \mbox{\rm for all}\; u\geq t\}\,.$$
We recall (see \cite{LPP2}) that
under the assumptions of the theorem,
there exists  $\GW$-a.s. an infinite
sequence of regeneration times $\{\tau_i\}_{i\geq 1}$, and 
the sequence $\{(|X_{\tau_{i+1}}|-|X_{\tau_i}|),
(\tau_{i+1}-\tau_i)\}_{i\geq 1}$ is i.i.d. under the $\GW$ measure,
and the variables
$|X_{\tau_2}|-|X_{\tau_1}|$ and $|X_{\tau_1}|$ possess exponential
moments
(see
\cite[Lemma 4.2]{DGPZ} for the last fact).  
A key to the proof of an annealed invariance principle
is the following
\begin{proposition}
	\label{prop-reg}
	When $\lambda<m$, it holds that $E_{\GW}((\tau_2-\tau_1)^k)<\infty$
	for all integer $k$.
\end{proposition}
{\bf Proof of Proposition \ref{prop-reg}:} 
By coupling with a biased (away from $0$)
simple random walk on $\mathbb{Z}_+$, the claim
is trivial if $\lambda<1$. The case $\lambda=1$ is covered
in \cite[Theorem 2]{piau}. We thus consider 
in the sequel only $\lambda\in (1,m)$.
Let
$T_o=\inf\{t>0: X_t=o\}$ denote the first 
return
time to the root and $T_n=\min\{t>0: |X_t|=n\}$ denote the 
hitting time of level $n$.
Let $o'\in D_1$ be an arbitrary offspring of the root.
By \cite[(4.25)]{DGPZ}, the law of $\tau_2-\tau_1$ under $\GW$ is
identical to the law of $\tau_1$ for the walk started at $v$, under the
measure $\GW^v(\cdot|T_o=\infty)$.
Therefore,
	$$E_{\GW}^o((\tau_2-\tau_1)^k)=
	E_{\GW}^{o'}(\tau_1^k\,| T_o=\infty)
	=   \frac{E_{\GW}^{o'}(\tau_1^k\,;\, T_o=\infty)}
{P^{o}_{\GW}(T_o=\infty)}
$$
where in the last equality we used that 
$P^{o}_{\GW}(T_o=\infty)=
P^{o'}_{\GW}(T_o=\infty)$.
Thus, with $c$ denoting a deterministic constant whose value
may change from line to line,
\begin{eqnarray*}
	E_{\GW}^o((\tau_2-\tau_1)^k)
	&\leq & c \sum_{n=1}^\infty  
	E_{\GW}^{o'}(\tau_1^k\,;\, |X_{\tau_1}|=n, T_o=\infty) \nonumber \\
	&= & c \sum_{n=1}^\infty  
	E_{\GW}^{o'}(T_n^k\,;\, |X_{\tau_1}|=n, T_o=\infty)  \\
	&\leq & c \sum_{n=1}^\infty  
	E_{\GW}^{o'}(T_n^{2k}\,;\,  T_o=\infty)^{1/2} 
	P_{\GW}^{o'}(|X_{\tau_1}|=n)^{1/2}\nonumber \\
	&\leq & 
	c \sum_{n=1}^\infty e^{-n/c} 
	E_{\GW}^o(T_n^{2k}\,;\,  T_o=\infty)^{1/2}\,,
\end{eqnarray*}
where the last inequality is due to the above mentioned exponential moments
on $|X_{\tau_1}|$. Therefore,
\begin{equation}
	\label{eq-190506a}
	E_{\GW}^o((\tau_2-\tau_1)^k)
	\leq  
	c \sum_{n=1}^\infty e^{-n/c}n^{10k} 
	\left(\sum_{j=0}^\infty (j+1)^{2k}
	P_{\GW}^o(T_n>jn^{10}\,;\,  T_o=\infty)\right)^{1/2}\,.
\end{equation}
We proceed by estimating the latter probability. For $j\geq 1$,
let
$${\cal A}_{1,j,n}=\{\mbox{\rm there exists a $t\leq jn^{10}$ such that}\,
d_{X_t}\geq (\log jn^{10})^2\}\,.$$
Note that by the assumption $\sum \beta^k p_k<\infty$ for some $\beta>1$, there
exists a constant $c$ such that for all $j$ and all $n$ large,
\begin{equation}
	\label{eq-190506b}
	P_{\GW}^o( {\cal A}_{1,j,n})\leq e^{-c (\log (jn^{10}))^2}\leq
	e^{-c(\log n^{10})^2-c (\log j)^2}\,,
\end{equation}
We next recall that $t$ is a fresh time for the random walk if
$X_s\neq X_t$ for all $s<t$. Let $N_{j,n}:=
|\{t\leq jn^{10}: \, \mbox{\rm $t$ is a fresh time}\}|\,$
(i.e., $N_{j,n}$ is the number of distinct vertices
visited by the walk up to time $jn^{10}$).
Set
$${\cal A}_{2,j,n}=
\{N_{j,n}< \sqrt{jn^{10}}\}\cap \{T_o=\infty\}\,.$$
Note that on the event 
${\cal A}_{2,j,n}\cap
{\cal A}_{1,j,n}^c$
there  is a time $t\leq
jn^{10}$ and a vertex $v$ with $d_v\leq (\log (jn^{10}))^2$
such that $X_t=v$ and $v$ is subsequently visited $\sqrt{jn^{10}}$ times 
with no visit at the root. Considering the trace of the walk on the ray 
connecting $v$ and $o$, and conditioning on $X_t=v$,
the last event has a probability bounded uniformly 
(in $t,v$) by $(1-c/(\log (jn^{10}))^2)^{\sqrt{jn^{10}}}$, since
$\lambda>1$. Hence,
for all $n$ large, using (\ref{eq-190506b}),
\begin{eqnarray}
	\label{eq-190506c}
	P_{\GW}^o( {\cal A}_{2,j,n})&\leq& e^{-c (\log (jn^{10}))^2}
	+
	jn^{10}
	\left(1-\frac{c}{(\log (jn^{10}))^2}\right)^{\sqrt{jn^{10}}}
	\nonumber \\
	&\leq &
	e^{-c(\log n^{10})^2-c (\log j)^2}+ jn^{10} e^{-(jn^{10})^{1/4}}
	\,.
\end{eqnarray}
The event 
${\cal A}_{2,j,n}^c\cap \{T_o=\infty\}$
entails the existence of at least $j^{1/2} n^3$ fresh  times
which are at distance at least $n^2$ from each other. Letting 
$t_1=\min\{t>0: t \, \mbox{\rm is a fresh time}\}$ and $$t_{i}=
\min\{t>t_{i-1}+n^2: t \, \mbox{\rm is a fresh time }\}\,,$$ we observe that
if $|X_{t_i}|<n$ then 
$P_{\GW}^{X_{t_i}}(T_{n}<n^2|{\cal F}_{t_i})>c>0$
(since from each fresh time, the walk has under the $\GW$ measure a strictly
positive probability to escape with positive speed without backtracking to the  
fresh point).  Thus,
\begin{equation}
	\label{eq-190506d}
	P_{\GW}^o(T_n>jn^{10}\,,\,  T_o=\infty, {\cal A}_{2,j,n}^c)
	\leq (1-c)^{j^{1/2} n^3}\,.
\end{equation}
Combining (\ref{eq-190506c}) and (\ref{eq-190506d}), we conclude that
$$	\sum_{j=0}^\infty (j+1)^{2k}
	P_{\GW}^o(T_n>jn^{10}\,,\,  T_o=\infty) \leq c\,.$$
Substituting in 
(\ref{eq-190506a}), the lemma follows.\qed

A standard consequence of Proposition
\ref{prop-reg} and the regeneration structure
(see e.g. \cite[Theorem 4.1]{sznitman},\cite[Theorem 3.5.24]{zeitouni}) 
is the following:
\begin{corollary}
	\label{cor-190506g}
	There exists a constant $\sigma^2$ such that,
	under the annealed measure $\GW$, the process
	$\{(|X_{\lfloor nt\rfloor}|-n\v t)/\sqrt{\sigma^2 n}\}_{t\geq 0}$ 
converges in distribution to a Brownian motion.
\end{corollary}

\noindent
{\bf Proof of Theorem \ref{theo-GWtransient}:}
Our argument is based on the technique introduced by 
Bolthausen and Sznitman in \cite{BS}, as developed in
\cite{BSZ}. Let $B^n_t=B^n_t(|X|_\cdot)=
(|X_{\lfloor nt \rfloor}|-nt \v)/\sqrt{n}$,
and let $\BBB^n_t(|X|_\cdot)$ denote the polygonal interpolation 
of $(k/n)\to B^n_{k/n}$. 
Consider the space 
$\CC_T$ of continuous functions 
on $[0,T]$, endowed with the distance
$d_T(u,u')=\sup_{t\leq T}|u(t)-u'(t)|\wedge 1$.
By \cite[Lemma 4.1]{BS},
Theorem \ref{theo-GWtransient} will follow
from Corollary \ref{cor-190506g}
once we show that for 
all bounded by $1$ Lipschitz function $F$ on
$\CC_T$ with Lipschitz constant $1$, and $b\in (1,2]$,
\begin{equation}
\label{eq-200805b}
\sum_k 
\mbox{\rm var}_{\GW}\left(E^o_{\cal T}
[F(\BBB^{\lfloor b^k
\rfloor})]\right)<\infty\,.
\end{equation}
In the sequel, fix $b$ and $F$ as above. For the same tree $\cal T$, let
$X^1_\cdot$ and $X^2_\cdot$ be independent $\lambda$-biased random walks on
$\cal T$, and set $\BBB[i,k]_t=\BBB^{\lfloor b^k\rfloor}_t(|X^i|_\cdot)$
and $\BBB[i,k,s]_t=\BBB^{\lfloor b^k\rfloor}_t(|X^i|_{\cdot+s}-|X^i|_s)$,
$i=1,2$.
Set
$$\tau^{i,k}=\min\{t>\lfloor b^{k/4}\rfloor:\,
\mbox{\rm $t$ is a regeneration
time for $X^i$}\}\, $$
$${\cal A}_k^1:=
\{\{X^1_s, s\leq \tau^{1,k}\}\cap X_{\tau^{2,k}}^2=\emptyset\},\quad
{\cal A}_k^2:=\{\{X^2_s, s\leq \tau^{2,k}\}\cap X_{\tau^{1,k}}^1=\emptyset\},$$
$${\cal A}_k={\cal A}_k^1\cap {\cal A}_k^2\,,$$
$$ {\cal B}_k^i:= \{ \tau^{i,k}\leq 
b^{k/3}\}. $$
Note that 
on the event ${\cal A}_k^1$, the paths 
$\{X^1_s, s\geq \tau^{1,k}\}$ and $\{X^2_s, s\geq \tau^{2,k}\}$ 
can intersect
only if $X^2_{\tau^{2,k}}$ is a descendant of $X^1_{\tau^{1,k}}$.
Applying the same reasoning for 
the symmetric event ${\cal A}_k^2$, we conclude
that on the event ${\cal A}_k$, these two paths do not intersect.

By construction, for any path $X_\cdot$ on ${\cal T}$,
the path $\BBB^{\lfloor b^k \rfloor}(|X|_\cdot)$ 
is Lipschitz with
Lipschitz constant bounded by $b^{k/2}$.
Hence,
since
$$\max_t|\BBB[i,k]_t-\BBB[i,k,\tau^{i,k}]_t|\leq \frac{\tau^{i,k}}{b^{k/2}}\,$$
and
using the fact that $F$ is a 
Lipschitz function with Lipschitz constant $1$, 
we have that
on the event ${\cal B}_k^i$, 
$|F(\BBB[i,k])-
F(\BBB[i,k,\tau^{i,k}])|\leq b^{k/3}/b^{k/2} 
$, and thus, since $|F|\leq 1$,
\begin{eqnarray*}
&&\mbox{\rm var}_{\GW}\left(E^o_{\cal T}
[F(\BBB^{\lfloor b^k
\rfloor})]\right)\\
&=& 
E_{\GW} [F(\BBB[1,k]) F(\BBB[2,k])]
-
E_{\GW} [F(\BBB[1,k])]E_{\GW} [F(\BBB[2,k])]\\
&\leq &
4P_{\GW}( ({\cal B}_k^1)^c)+ 4b^{k/3-k/2}
+E_{\GW} [F(\BBB[1,k,\tau^{1,k}]) F(\BBB[2,k,\tau^{2,k}])]
\\&&\quad \quad -
E_{\GW} [F(\BBB[1,k,\tau^{1,k}])]E_{\GW} 
[F(\BBB[2,k,\tau^{2,k}])]\,.\end{eqnarray*}
Conditioning on the event ${\cal A}_k$ and using again
that $|F|\leq 1$, we get
\begin{eqnarray*}
\mbox{\rm var}_{\GW}\left(E^o_{\cal T}
[F(\BBB^{\lfloor b^k
\rfloor})]\right) 
&\leq &
4P_{\GW}( ({\cal B}_k^1)^c)+4P_{\GW}({\cal A}_k^c)+4b^{k/3-k/2}\\&&\quad
+
E_{\GW} [F(\BBB[1,k,\tau^{1,k}]) F(\BBB[2,k,\tau^{2,k}])| {\cal A}_k]
\\&&\quad -
E_{\GW} [F(\BBB[1,k,\tau^{1,k}])|{\cal A}_k]E_{\GW} 
[F(\BBB[2,k,\tau^{2,k}])|{\cal A}_k]
\,.
\end{eqnarray*}
Conditioned on the event ${\cal A}_k$, the paths
$\BBB[1,k,\tau^{1,k}]$ and $\BBB[2,k,\tau^{2,k}]$ 
are independent under
the $\GW$ measure. Therefore, we conclude that
\begin{equation}
	\label{eq-200506a}
	\mbox{\rm var}_{\GW}\left(E^o_{\cal T}
[F(\BBB^{\lfloor b^k
\rfloor})]\right)\leq  4(P_{\GW}({\cal A}_k^c)+
P_{\GW}( ({\cal B}_k^1)^c)+b^{k/3-k/2})\,.
\end{equation}
Let $\tau_i^j$ denote the successive regeneration times for
$X_\cdot^j$, $j=1,2$. 
The event 
$\{({\cal B}_k^1)^c\}\cap 
\{\tau_1^1\leq b^{k/4}\}$ implies that at least one of the first
$b^{k/4}$ inter-regeneration times $\tau_{i+1}^1-\tau_i^1$ is larger than
$b^{k/3}$. Therefore,
\begin{eqnarray*}
	P_{\GW}( ({\cal B}_k^1)^c)&\leq&
	P_{\GW}(\tau_1^1>b^{k/4})+
	P_{\GW}(\mbox{\rm one of $(\tau_{i+1}^1-\tau_i^1)_{i=1}^{b^{k/4}}$
 is larger than 
$b^{k/3}$})\\
&\leq&  P_{\GW}(\tau_1^1>b^{k/4})+
b^{k/4}P_{\GW}(\tau_2^1-\tau_1^1>b^{k/3})\\
&\leq& 
P_{\GW}(\tau_1^1>b^{k/4})+
b^{k/4-k/3} E_{\GW}(\tau_2^1-\tau_1^1)\\&\leq&
P_{\GW}(\tau_1^1>b^{k/4})+ c b^{k/4-k/3}
\,. 
\end{eqnarray*}
where Markov's inequality was used in the third step.
Let $T_\ell=\min\{t>0: |X_t|=\ell\}$.
Let $Y_t$ be a nearest neighbor random walk on 
$\mathbb{Z}_+$ with $P(Y_{t+1}=Y_t-1|Y_t)=\lambda/(\lambda+1)$ whenever
$Y_t\neq 0$. $Y_\cdot$ and $X_\cdot$ can be constructed on the same
probability space, such that $T_\ell\leq \min\{t>0: Y_t=\ell\}=:T^Y_\ell$ 
for 
all $\ell$. On the other hand, using the Markov property,
for any constant $c$ and all $\ell$ large,
$$P(T^Y_\ell>e^{c\ell})\leq 
 \left(1-\left(\frac{1}{1+\lambda}\right)^\ell\right)^{e^{c\ell}/\ell}
$$
In particular, 
there exists a $c_1=c_1(\lambda)>0$ such that
$P_\GW(T_\ell>e^{c_1\ell})\leq e^{-\ell/c_1}$ (better bounds are available
but not needed). Thus, for some deterministic constants 
$c_i=c_i(\lambda,b)>0$, $i\geq 2$, 
and all
$k$ large,
\begin{eqnarray}
	\label{eq-200506g}
	P_{\GW}(\tau_1>b^{k/4}/2)& \leq &
\!\!\!\!
P_{\GW}(|X_{\tau_1}|> c_2 k)+
P_{\GW}(\tau_1>b^{k/4}/2,
|X_{\tau_1}|\leq c_2 k)\nonumber \\
&
\leq&
\!\!\!\!
P_{\GW}(|X_{\tau_1}|> c_2 k)+
P_{\GW}(T_{c_2 k}> b^{k/4}/2)
\leq e^{-c_3 k},
\end{eqnarray}
where we have used the above mentioned fact that 
$|X_{\tau_1}|$ possesses exponential moments.
We conclude that with  $c_4\leq c_3$,
$$P_{\GW}( ({\cal B}_k^1)^c)\leq
 b^{-c_4 k}\,.$$

It remains to estimate $P_{\GW}({\cal A}_k^c)
\leq 2P_{\GW}( ({\cal A}_k^1)^c)$.
Let $${\cal C}_{k,i}':=\{\tau_1^i<b^{k/4}/2\},\,
{\cal C}_{k,i}'':=\{
\tau_{\lfloor b^{k/8}\rfloor}^i<b^{k/4}\},\,
{\cal C}_k:={\cal C}_{k,1}'\cap {\cal C}_{k,2}'\cap {\cal C}_{k,1}''\cap 
{\cal C}_{k,2}''\,.$$
Using (\ref{eq-200506g}), it follows that
$P_\GW({\cal C}_{k,1}'^c)\leq 
	b^{-c_3 k}$.
	On the other hand, the event ${\cal C}_{k,1}'\cap ({\cal C}_{k,1}'')^c$ implies
	that the sum of the difference $\tau^1_{i+1}-\tau^1_i$, $i=1,\ldots,
	\lfloor b^{k/8}\rfloor$, is larger than $b^{k/4}/2$, and hence, by
	Markov's inequality,
	$$P_{\GW}(
	{\cal C}_{k,1}'\cap ({\cal C}_{k,1}'')^c)\leq 
	2b^{k/8}\frac{E_{\GW}(\tau_2^1-\tau_1^1)}
{b^{k/4}}\leq b^{-c_5 k}\,,$$
for some deterministic constant $c_5<c_4$. Since the same estimates
are valid also for $C_{k,2}'$ and $C_{k,2}''$ replacing
$C_{k,1}'$ and $C_{k,1}''$,
it follows that
\begin{equation}
	\label{eq-200506h}
	P_{\GW}({\cal C}_k^c)
	\leq 
4b^{-c_5 k}\,.
\end{equation}
On the other hand, let ${\cal Z}^i$ denote the collection of vertices in 
$D_{\lfloor b^{k/8}\rfloor}$ hit by $X^i_\cdot$. On ${\cal C}_k$ there
are at most $b^{k/4}$ vertices in ${\cal Z}^1$. 
The event $({\cal A}_k^1)^c \cap {\cal C}_k$ implies that
the path $X^2$ intersected the path $X^1$ at a distance at least 
$\lfloor b^{k/8}\rfloor$ from the root, and this has to happen
before time $\tau^2_{\lfloor b^{k/8}\rfloor}$, i.e. before  time
$b^{k/4}$, for otherwise ${\cal Z}^1\cap {\cal Z}^2=\emptyset$.
Therefore,
\begin{eqnarray}
	\label{eq-200506j}
	P_{\GW}( ({\cal A}_k^1)^c\cap {\cal C}_k)
	&\leq &E_{\GW}P_{\cal T}^o( X^2_\cdot \, \mbox{\rm visits ${\cal Z}^1$
	before time $b^{k/4}$})
 \\ &	\leq& b^{k/4}
	E_{\GW}\max_{v\in D_{\lfloor b^{k/8}\rfloor}}
P_{\cal T}^o( 
	X^2_\cdot \, \mbox{\rm visits $v$
	before time $b^{k/4}$})\,.
\nonumber
\end{eqnarray}
When $\lambda>1$, there exists a constant $c_6<c_5$ such that
uniformly in  $v\in  
D_{\lfloor b^{k/8}\rfloor}$,
$$P_{\cal T}^o( 
	X^2_\cdot \, \mbox{\rm visits $v$
	before time $b^{k/4}$}) \leq b^{k/4}e^{-c_6b^{k/8}}\,.$$ 
	On the other hand, even when $p_1>0$, Lemma 2.2 of \cite{DGPZ}
	shows that there exists a $\beta>0$ such that
	with $M_v=|\{w \, \mbox{\rm is an ancestor of $v$}: d_w\geq 2\}|$,
	it holds that
	$$\limsup_{\ell\to\infty}P_\GW(\min_{v\in 
	D_\ell} M_v/\ell<\beta)<0\,.$$ 
	It immediately follows, reducing $c_6$ if necessary,
	that when $\lambda\leq 1$, for all $k$ large,
$$		E_{\GW}\max_{v\in D_{\lfloor b^{k/8}\rfloor}}
P_{\cal T}^o( 
X^2_\cdot \, \mbox{\rm visits $v$ ever}
	)\leq  e^{-c_6 b^{k/8}}\,.$$
	Substituting in (\ref{eq-200506j}), we conclude that
	whenever $\lambda< m$,
	$$ 
	P_{\GW}( {\cal A}_k^c\cap {\cal C}_k)\leq 2
	P_{\GW}( ({\cal A}_k^1)^c\cap {\cal C}_k)\leq 
	 e^{-c_7 b^{k/8}}\,.$$
	 Together with (\ref{eq-200506h}),
	  (\ref{eq-200506g}), and
	  (\ref{eq-200506a}),  we conclude that
(\ref{eq-200805b}) holds  and thus conclude the proof of 
Theorem \ref{theo-GWtransient}. \qed

\noindent
{\bf Acknowledgment} We thank Nina Gantert for asking the question 
that led to this work, and for many useful discussions. We thank
Nathan Levy for  a careful reading of several  earlier versions of this 
paper.

\end{document}